\documentclass[11pt,a4paper]{article}
\usepackage{epsfig}
\usepackage{amsmath}
\usepackage{amsfonts}
\usepackage{amssymb}
\usepackage{amsthm}

\newcommand{\old}[1]{}
\newcommand{\mn}{\medskip\noindent}

\newcommand{\ga}{\alpha}
\newcommand{\gb}{\beta}
\newcommand{\gd}{\delta}

\newcommand{\bO}{{\mathcal O}}

\newcommand{\bn}{\bigskip\noindent}
\newcommand{\lm}{\ell_m}

\newtheorem{theorem}{Theorem}[section]

\newtheorem{lemma}[theorem]{Lemma}
\newtheorem{corollary}[theorem]{Corollary}

\begin{document}

\title{Finding  Weighted Graphs  by   Combinatorial Search }

\author{Jeong Han Kim\footnote{Department of Mathematics, Yonsei University, 
Seoul, 120-749 Korea (e-mail: jehkim@yonsei.ac.kr).}
 }

\date{}

\maketitle

\newcommand{\beq}{\begin{equation}}
  \newcommand{\enq}{\end{equation}}
\newcommand{\raf}[1]{(\ref{#1})}
\newcommand{\pr}{\Pr}
\newcommand{\el}{{\! \ell}}
\newcommand{\eps}{\varepsilon}
\newcommand{\gs}{\sigma}

\begin{abstract}

We consider the problem of finding edges of  a hidden weighted   graph  using a certain type of queries. Let $G$ be a weighted graph with $n$ vertices.
In the most general setting,  the $n$ vertices are known and no other information about $G$  is given. The problem is finding all edges of $G$ and their weights using  additive queries, where, for an additive query, one chooses a set of vertices and asks the sum of the weights of edges with both ends in the set. This model has been  extensively used in bioinformatics including genom sequencing. Extending recent results of Bshouty and Mazzawi \cite{BM11_TCS}, and Choi and Kim \cite{CK12}, we present a polynomial time randomized algorithm to find the hidden weighted graph $G$ when the number of edges in $G$ is known to be at most $m\geq 2$ and the weight $w(e)$ of each edge $e$ satisfies
$\ga \leq |w(e)|\leq \gb$ for fixed constants $\ga, \gb>0$. The query complexity of the algorithm is $O(\frac{m \log n }{\log m})$, which is optimal up to a constant factor.

The algorithm heavily relies on a well-known combinatorial search problem, which may be of independent interest. Suppose that there are $n$ identical looking coins and  some of them are counterfeit. The weights of all authentic coins are the same and known a priori. The weights of counterfeit coins vary but different from the weight of an authentic coin. Without loss of generality, we may assume the weight of authentic coins is $0$.
The problem is to find all counterfeit coins by weighing sets of coins on a spring scale. We introduce a polynomial time randomized algorithm to find all counterfeit coins when the number of them is know to be at most $m\geq 2$ and the weight $w(c)$ of each counterfeit coin $c$ satisfies
$\ga \leq |w(c)|\leq \gb$ for fixed constants $\ga, \gb$.
The query complexity of the algorithm is $O(\frac{m \log n }{\log m})$, which is optimal up to a constant factor. The algorithm uses, in part, random walks.

\mn


\end{abstract}

\section{Introduction}

\subsection{Graph Finding Problem}
We consider the problem of finding edges of  a hidden weighted graph  using a certain type of queries. Let $G$ be a weighted graph with $n$ vertices.
In the most general setting,  the $n$ vertices are known and no other information about $G$  is given. The problem is finding all edges of $G$ and their weights using  queries. Three types of queries have been extensively studied:

\mn
{\bf Detection query:} One chooses a set of vertices and asks  if there is an edge with both ends in the set. This type of queries has applications to genom sequencing and has been studied in \cite{AA05,ABKRS04,AC04,AC06,GK97,GK98}.

\mn
{\bf Additive query:} One chooses a set of vertices and asks  the sum of weights of edges with both ends in the set. This model has been  extensively used in bioinformatics including genom sequencing, and  studied in \cite{AC04,BGK05,BM10_STACS,BM10_MFCS,BM11_SODA,BM11_TCS,CK10_AI, Grebinski98,GK98,GK00,Mazzawi10,RS07}.

\mn
{\bf Shortest path query:} One choose a pair of vertices and asks the length of the shortest path between the two vertices. This query arises in  the canonical model
of the evolutionary tree literature \cite{Hein89,KZZ03,RS07a}.

\mn
(Our lists of references are far from being exhaustive.)

\mn

In this paper, we focus on the additive queries. The graph finding problem with additive queries is partly motivated by the shotgun sequencing \cite{BAFK01,GK98}, one of the most popular methods for DNA sequencing. In the shotgun sequencing, one needs
to put back separately decoded short fragments of a given genome sequence
into the same order as in the original sequence.
Combined with a biotech method called the multiplex PCR \cite{TRKKS99},
the process is reduced to the problem of finding a hidden graph using  additive queries.
The additive queries are also  used in the problem of finding the Fourier coefficients of pseudo-Boolean functions, which play crucial roles in evolutionary computation, artificial intelligence, and population genetics \cite{CK10_AI,CJK11_JCSS,CK11_AI}.

In the rest of this paper, we say queries for  additive queries and
all logarithms are in base 2, unless otherwise specified.
For unweighted graphs,
Grebinski and
Kucherov presented a few results. For arbitrary graphs on $n$ vertices, they have shown that $O(\frac{n^2}{\log n})$ queries are enough \cite{GK00}. If the hidden graph is known to be a Hamiltonian path or cycle, then $O(n)$ queries are suffice \cite{GK98}.
More generally, if the maximum degree of the hidden graph is known to be at most $d$, then  the graph may be found using $O(dn)$ queries \cite{GK00}. Grebinski \cite{Grebinski98} has shown that the same bound $O(dn)$ holds for d-degenerate graphs.

When the hidden graph has at most $m\geq 2$ edges and $m$ is known,
some bounds close to the optimal bound were shown \cite{AC04, RS07} and Choi and Kim \cite{CK10_AI} proved a $O(\frac{m\log (n^2/m)}{\log m})$ bound that is optimal (up to a constant factor). The  randomized algorithm presented there uses non-adaptive queries but it is not a polynomial time algorithm, where  queries are non-adaptive if each query is independent of answers to the previous queries.
 Recently, Mazzawi
\cite{Mazzawi10} constructed a polynomial time algorithm with optimal query complexity. The algorithm is deterministic and  uses  adaptive queries. She also extended the algorithm to find weighted graphs with positive integer weights.

For weighted graphs,
Choi and Kim \cite{CK10_AI} proved a non-adaptive $O(\frac{m\log n}{\log m})$ query bound, provided that $m$ is at least a polylog of $n$ and the absolute values of  all weights are between $n^{-a}$ and $n^b$ for constants $a,b>0$. Bshouty and Mazzawi \cite{BM11_TCS} showed the same bound without the extra conditions. However, it is unlikely that one may able to develop a polynomial time algorithm from those results. In other words,  substantially new ideas seem to be needed to design an algorithm that is useful in practical sense. A significant result toward this direction has been  shown by Bshouty and Mazzawi \cite{BM10_MFCS}: For weighted graphs with positive real weights, they presented a deterministic polynomial time algorithm that uses an almost optimal number of (adaptive) queries,
$O(\frac{m\log n}{\log m} + m\log\log m)$.
Note that the extra $m\log\log m$ term is larger than the optimal query bound by a $\log\log n$ factor when $\log m =\Omega (\log n)$.

 To obtain the optimal query complexity $O(\frac{m\log n}{\log m})$, Choi and Kim \cite{CK12} have recently introduced  a randomized polynomial time algorithm that finds the hidden weighted graph with positive real weights.
Another randomized polynomial time algorithm
they introduced uses $O(\frac{m\log n}{\log m})$ queries to find  the hidden weighted graph with bounded integer weights.

In this paper, we present a randomized polynomial time algorithm
that  works for a quite general class of weighted graphs.
Using the optimal number of queries up to constant factor, the algorithm finds the hidden weighted graph provided that  the weight $w(e)$ of each edge $e$ in the graph satisfies $\ga\leq |w(e)| \leq \gb$ for positive constants $\ga$ and $\gb$. The theorem we will prove is slightly more general in the sense that $\ga, \gb$  are not necessarily constants.

\begin{theorem}\label{gfpm}
Let $n,m$ be positive integers with $n^2\geq m\geq 2$ and let $\ga, \gb>0$ be  positive real numbers (not necessarily constants) with $2\ga<\gb$.
Suppose  a weighted graph $G$ with $n$ vertices and at most $m$ edges is given.   If the weight $w(e)$ of each edge in $G$ satisfies $\ga \leq |w(e)|\leq \gb$, then there is a randomized polynomial time algorithm that asks $O(\frac{m \log (\gb/\ga) \log n}{\log m})$ queries to
 find all edges with probability $1-O(1/{m^{0.02}})$.
\end{theorem}

\old{It is not clear that the lower and upper bounds for the absolute values of edge weights are really necessary. It seems to be not extremely
surprising even if they were, though.}

Our proof of the theorem heavily relies on a well-known combinatorial search problem. Suppose there are $n$ identical looking coins and  some of them are counterfeit. The weights of all authentic coins are the same and known a priori. The weights of counterfeit coins vary but different from the weight of an authentic coin.
The problem is to find all counterfeit coins by weighing sets of coins on a spring scale.
Note that  weighing sets of coins on a spring scale may be regarded as additive queries. This problem is also equivalent to the  graph finding problem when the graphs are restricted to stars $K_{1,m} $ with  known center. The coin weighing problem has been extensively studied. We survey its colorful history and add one more algorithm  finding all counterfeit coins when the weights of each counterfeit coin satisfies  properties similar to those described in the above theorem.

\subsection{Coin Weighing Problem}

Suppose  there are $n$ identically looking coins, some of them are counterfeit. The weights of all authentic coins are the same and known a priori, while the weights of  counterfeit coins are unknown but different from  the weight of an  authentic coin. Without loss of generality, it may be assumed that the weights of authentic coins are $0$ and the weights of counterfeit coins belong to a set  of non-zero real numbers. We want to find all counterfeit coins by weighing sets of coins on a spring scale, which we call additive queries or simply queries.

After the coin weighing problem was introduced by Fine \cite{Fine60} and Shapiro \cite{Shapiro60}, a number of results have been published,
mainly focusing on the case that the weights of counterfeit coins are the same \cite{Cantor64,CM66,ER63,Lindstrom64,Lindstrom65,Lindstrom71, Moser70,SS63}: Summarizing some of them briefly, Erd\H{o}s and
R\'{e}nyi \cite{ER63},  in 1963, proved that $\frac{(\log 9 +o(1))n}{\log n}$ queries are enough and $\frac{(2+o(1))n}{\log n}$ queries are required.
(See  \cite{LV94} for another proof of the lower bound.) The upper bound was improved to match the lower bound by Cantor and Mills \cite{CM66}, and Lindstr\"{o}m
\cite{Lindstrom65}. Using the M\"obius function,  Lindstr\"{o}m \cite{Lindstrom71,Lindstrom75} explicitly constructed a query matrix  that asks $\frac{(2+o(1))n}{\log n}$ queries.
The case that  the number $m$ of counterfeit coins is also known has been   extensively studied too \cite{Capetanakis79b,Capetanakis79a,
DH93,GK00,Lindstrom75,M81,TM78,UTW00}. Recently,  Bshouty \cite{Bshouty09} proposed the first polynomial time  algorithm
 that uses $\frac{(1+o(1))2m\log \frac{n}{m}}{\log m}$ adaptive queries. The query complexity is  optimal up to $o(1)$ term.
 \old{ When the weights of all coins are non-negative integers and their sum  is known (instead of the number $m$ of counterfeit coins), then he was able to modify  the algorithm to  find all counterfeit coins with an optimal query complexity  up to a constant factor.}

Results for the general case, in which the weights of counterfeit coins are not the same, have been obtained only recently. As the results were applied to the (weighted) graph finding problem, our summary is almost the same as in the previous subsection. When the weights of the counterfeit coins can be any (not necessarily positive) real numbers,
 Choi and Kim \cite{CK10_AI} proposed an algorithm with a non-adaptive $O(\frac{m\log n}{\log m})$ query bound, under the mild conditions on $m$ and the weights, i.e.,   $m=\Omega({\rm polylog} n)$ and the absolute values of  all weights are between $n^{-a}$ and $n^b$ for constants $a,b>0$. Bshouty and Mazzawi \cite{BM11_TCS} showed the same bound without the extra conditions. Though the query complexities of both algorithms are optimal, the time complexities of them are far from being polynomial. Concerning polynomial time algorithms, Bshouty and Mazzawi \cite{BM10_MFCS} presented a deterministic polynomial time algorithm that uses a near optimal number of (adaptive) queries,
$O(\frac{m\log n}{\log m} + m\log\log m)$, assuming the weights of all counterfeit coins are positive real numbers. They first constructed    a search matrix using Fourier representations, and took the divide and conquer approach to guess  the  sums of the weights of coins. The search matrix played key roles when the sums of the weights were guessed. The processes for checking and correction follow after guessing.

As mentioned before, the extra $m\log\log m$ term is larger than the optimal bound by a  $\log \log n$ factor when $\log m=\Omega (\log n)$. Choi and Kim \cite{CK12} presented a polynomial time randomized  algorithm to remove the $m\log\log m$ term in the query complexity. Another polynomial time randomized  algorithm may be applied to achieve the optimal query complexity, when the weights of counterfeit coins are bounded integers
  in absolute values. The key idea is constructing random sets of coins that are useful to control the number of checking and correction processes used by Bshouty and Mazzawi \cite{BM10_MFCS}.
  Once the number of checking and correction processes is substantially reduced, less queries are needed.

A randomized algorithm is presented in this paper to achieve the optimal query complexity when the weights of counterfeit coins are any real numbers bounded from below and from above in absolute values.
The theorem we will prove is slightly more general in the sense that some exceptions for the weight condition  are allowed.

\begin{theorem}\label{cw} Let $n,m$ be positive integers with $n\geq m\geq 2$ and let $\ga, \gb, \eps >0$ be  positive real numbers (not necessarily constants) with $2\ga<\gb, \eps <1/2$.
Suppose  $n$ coins are given and there are at most $m$ counterfeit coins among them. The weights of authentic coins are $0$ while the wights of counterfeit coins vary but they are non-zero.
 If the weights $w(c)$ of all but $\eps m$ counterfeit coins $c$ satisfy $\ga \leq |w(c)|\leq \gb$ and the weights $w(c)$ of the $\eps m$ counterfeit coins $c$ satisfies just $|w(c)| \leq \gb$, then there is a randomized polynomial time algorithm that asks $O(\frac{m \log (\gb/\ga) \log n}{\log m})$ queries and finds all but $m^{0.8}+2\eps m$ counterfeit coins, with probability $1-O(1/{m^{0.8}})$.
All the remaining counterfeit coins can be found using
  $O((m^{0.8}+2\eps m)\log n)$ additional queries, with probability
$1-e^{-\Omega(m^{0.8})}$.
\end{theorem}

\noindent
In the proof of Theorem \ref{cw}, we use the search matrix Bshouty and Mazzawi \cite{BM10_MFCS} developed after constructing random sets of coins as in Choi and Kim \cite{CK12}. Though the guessing processes are the same as in \cite{CK12}, the processes for checking and correction are newly developed using biased random walks.

One may easily verify if the coins declared  to be counterfeit by the algorithm in Theorem \ref{cw}  are actually counterfeit by directly weighing, using  $m$ additional queries. Running the algorithm $O(\mu)$ times with the verification at each time, the error probability may be made arbitrarily small.

\begin{corollary}\label{cwc}
Under the same hypotheses of Theorem \ref{cw} and any integer $\mu \geq 1$, there is a randomized polynomial algorithm that uses  $O(\frac{\mu m \log (\gb/\ga) \log n}{\log m})$ queries and finds all but $m^{0.8}+2\eps m$ counterfeit coins with probability $1-O(1/m^{\mu})$. All the remaining counterfeit coins can be found using
  $O((m^{0.8}+2\eps m)\log n)$ additional queries, with probability
$1-e^{-\Omega(m^{0.8} )}$.
\end{corollary}

\mn

After presenting the search matrix and two martingale  inequalities
in Section \ref{pre}, we prove Theorem \ref{cw} in Section \ref{scwp}. Section \ref{sgfp} is for the proof of Theorem \ref{gfpm}. The concluding remark will follow.

\section{Preliminaries}\label{pre}

As mentioned in the previous section,  Bshouty and Mazzawi \cite{BM10_MFCS} used Fourier representation of
certain functions to find a search matrix, i.e., a $0,1$ matrix that
is useful for coin weighing problems. We present properties of the
matrix in a slightly generalized form.

\newcommand{\gc}{\gamma}
\begin{lemma}\label{BM}  Let  $\gc, m$ be positive integers. Then,
for the smallest integer $t$ satisfying $t2^{t-1} \geq \gc m$, one
can  construct, in polynomial time, $2^t \times m$ $0,1$ matrix $S$
and $ 2^t \times 2^t$ matrix $T$ with the following property: For
each   $j=1,...,m$,  one may find, in polynomial time, a unique
positive integer $i_j \leq 2^{t}$ and a non-negative integer $k_j
\leq \lceil t/\gc \rceil -1 $  satisfying
$$(TS)_{i_{\! j} k} = 2^{-(k-j) \ga}
(TS)_{i_{\! j} j} ~\mbox{ for ~$j+1\leq k\leq j+k_j $},~~{\rm and}~~
(TS)_{i_{\! j} k}=0 ~\mbox{ for  ~$ k\geq j+k_j+1 $}, $$ where
$(TS)_{ij}$ is the $ij$ entry of $TS$.
\end{lemma}

Setting $ a_{jk} = \frac{(TS)_{i_{\! j} k}}{(TS)_{i_{\! j} j}}$, we
have the following corollary.

\begin{corollary}\label{sm}  Let $\gc, m$ are positive integers and $t$ be
the smallest integer  satisfying $t2^{t-1}\geq\gc m$. Then one can
find, in polynomial time, $2^t$ non-adaptive queries, real numbers
$a_{jk}$, and  a non-negative integer $k_j \leq \lceil t/\gc \rceil
- 1$, $ j=1,..., m, k=1,..., j-1$, satisfying the following
property: For disjoint sets $A_1,..., A_m$ of coins, the $2^t$
queries yield values $x_j$, in polynomial time,
satisfying
$$
  w(A_j) = x_j  -\sum_{k=1}^{j-1}  a_{j k} w(A_k)   - \sum_{k=1}^{k_j} \frac{
  w(A_{j+k}) }{2^{k\ga}},
  $$
  $j=1,..., m$, where $w(A)$ is the sum of weights of all coins in $A$.
  In particular, $\frac{(2+o(1))\gc
m}{\log (\gc m)}$ queries are enough to find $x_{j}$'s.
\end{corollary}

We will need the  Azuma-Hoeffding  martingale inequality too.
The following is from \cite{McDiarmid89}.

\begin{lemma}\label{mar}
Let $Z=(Z_1,\ldots,Z_t)$ be a family of independent random variables
with $Z_\ell$ taking values in a finite set $B_\ell$ for each $\ell$. Suppose
that the real-valued function $f$ defined on $\prod_\ell B_\ell$ satisfies
\begin{equation*}
|f(\mathbf{z})-f(\mathbf{z}')| \leq c_{_\ell}
\end{equation*}
whenever the vectors $\mathbf{z}$ and $\mathbf{z}'$ differ only in
the $\ell^{\rm th}$ coordinate. Then for any $\lambda \geq 0$,
\begin{equation*}
\Pr \left[ |f(Z) - \mathrm{E}[f(Z)]| \geq \lambda \right] \leq
2e^{-2\lambda^2 / \sum_{\ell} c_{_\ell}^2} .
\end{equation*}
\end{lemma}

For our purpose, a more general martingale inequality is needed. The
following version  appeared in \cite{Kim95}.
\begin{lemma}\label{gm}
Let $X=(Z_1,\ldots,Z_t)$ be independent identically distributed {\em
(}i.i.d.{\em )} Bernoulli random variables with probability $p$ {\em
(}i.e., $\Pr [Z_i = 1] = p$ and $\Pr [Z_i = 0] = 1-p$ for each
$i${\em )}. Suppose that the real-valued function $f$ defined on
$\{0,1\}^{t}$ satisfies
\begin{equation*}
|f(\mathbf{z})-f(\mathbf{z}')| \leq c_i
\end{equation*}
whenever the vectors $\mathbf{z}$ and $\mathbf{z}'$ differ only in
the $i^{\rm th}$ coordinate. Then for any $\lambda,\rho > 0$,
\begin{equation*}
\Pr \left[ |f(Z) - \mathrm{E}[f(Z)]| \geq \lambda \right] \leq 2
\exp \Big( -\rho \lambda + (\rho^{2}/2)p(1-p) \sum_{i=1}^{t} c_i^2
\exp(\rho c_i) \Big) .
\end{equation*}
\end{lemma}

\section{Coin Weighing Problem}\label{scwp}
\newcommand{\q}{\ell_q}
\newcommand{\lmq}{\ell_q}

Suppose  $n$ coins are given, some of which are counterfeit.
The weights of all authentic coins are the same and known a priori, while the weights of  counterfeit coins are unknown but different from  the weight of an  authentic coin. Without loss of generality, we may assume that the weights of authentic coins are $0$ and the weights of counterfeit coins belong to a set  of non-zero real numbers. We assume that the number of counterfeit coins is known to be at most $m$.

If $O(m\log n)$ queries are allowed to find counterfeit coins. One may use a randomized binary search:

\mn
{\bf Randomized Binary Search} Suppose a set $A$ of coins is given, and the number of coins is no more than $n$ and there are  at most $m\leq n$ counterfeit coins. Then select each coin with probability $1/2$, independently of all other coins. Then weigh the set $A'$ of selected coins. If the weight is non-zero, then find a counterfeit coin among the selected coins, using the deterministic binary search.

\mn
The deterministic binary search  is as follows. Divide $A'$ into two parts $A_{1}'$,   $A_{2}'$ with size difference at most $1$. If $w(A_{1}')\not=0$, then select $A_{1}'$. Otherwise, select $A_{2}'$.
Keep doing this for the selected set until  a counterfeit coin is found.

 \mn

 Provided that there is a counterfeit coin, it is not hard to see  that  the probability of the weight of $A'$  being  non-zero is at least $1/2$ and the deterministic binary search requires no more than $\lceil \log n\rceil$ queries. The number of queries required to find one counterfeit coin is at most $2+\lceil \log n\rceil$ in expectation. Thus, it is expected that $(\lceil \log n\rceil+2 +o(1))m$ queries are enough to find all counterfeit coins, with hight probability. Here, we show that $(\lceil \log n\rceil+3 )m$ queries are enough, with probability $1-e^{-\Omega(m)}$.

\begin{lemma}\label{rbs}  With probability $1-e^{-\Omega(m)}$, the randomized binary search   finds all counterfeit coins using $(\lceil \log n\rceil+3 )m$ queries.
\end{lemma}

\noindent
The proof of the lemma is presented in Appendix.

\mn

We first construct random sets of coins and then present the algorithm,  for which the time complexity is not optimized but it is clearly a polynomial time algorithm. Some explanation and analysis of the algorithm will follow after the algorithm is presented.
The construction of random sets is the same as in Choi and Kim \cite{CK12}.

\mn {\bf Constructing random sets of coins}: Let $A$ be a set of  $n$ or less coins. For an integer $q\geq 2$ and $\ell_q :=\lceil
\log q \rceil$,
we  construct random subsets $A_{i,j}$ of $A$,
  $i=0,1, ..., \lceil 3 \log n \rceil$,
  $j=1,..., 2^{\ell_q +i}  $.
For  $i=0$, we  assign  each coin in $A$ a uniform random number among $1,
..., 2^{\q}$, independently of all other coins. The set $A_{0,j}$
consists of all coins with assigned number $j$. Generally, for $i=1,...,\lceil 2\log q \rceil-1$,  once all
$A_{i-1,j}$, $j=1, ..., 2^{\lmq+i-1}$, are constructed, we may randomly
divide each set $A_{i-1,j}$ into two parts so that coins in $A_{i-1,j}$
are independently in the first part  with probability $1/2$. The
other coins in $A_{i-1,j}$ are to be in the second part. The set of
all coins in the first and second parts are denoted by $A_{i,2j-1}$
and $A_{i,2j}$, respectively. Or equivalently, after assigning
each coin mutually independent random numbers $r_{_0}, r_{_1}, ...,
r_{_{\lceil 2\log q \rceil-1}}$, independently of all other coins,
with
$$ \Pr [ r_{_{\!0}}= a ]= 2^{-\q},~ a=1, ..., 2^{\lmq}~~\mbox{and}~~~
\Pr [ r_{_{\!i}}= a ]=\frac{ 1}{2}, ~~a=0, 1,~~i=1,...,\lceil 2\log
q \rceil-1,
$$
we define  $A_{i,j}$ to be the set of all coins with assigned
numbers $r_{_{\!0}}, r_{_{\!1}}, ...., r_{_{\lceil 2\log q \rceil-1}}$
satisfying $j=1+(r_{_{\! 0}}-1)2^i +r_{_{\! 1}} 2^{i-1}+ \cdots +
r_{_{\! i}} $.

For $i  \geq  \lceil 2\log q \rceil$, $A_{i-1,j}$  may be deterministically
divided
 into two parts so that the first part consists of $\lceil
|A_{i-1,j}|/2 \rceil$ coins. As before, the first part is denoted by
$A_{i,2j-1}$, and  $A_{i,2j}=A_{i-1,j}\setminus A_{i,2j-1}$. This
construction is expected to stop when all $A_{i, j}$,
$j=1,...,2^{\lm+i}$, consist of one or no coin. As there are $n$ coins,  all $A_{i, j}$  consist of one or no coin within $\lceil \log n \rceil$ more rounds after
the random construction ends. For the sake of safeness, we stop the construct when $i=\lceil 3\log n \rceil
\geq \lceil 2\log q \rceil + \lceil \log n \rceil. $

\mn

The following lemma summarize properties of the random subsets $A_{i,j}$ that will be used for the analysis of the algorithm presented later. The proof is essentially in \cite{CK12} and it is presented in Appendix for the sake of completeness.

\begin{lemma}\label{random}
Suppose a set $A$ of   $n$ or less coins are given, and the number of counterfeit coins in $A$ is at most $q\geq 2$. If the weights $w(c)$ of all but at most $q/2$ counterfeit coins $c$ satisfy $|w(c)|\geq \ga$.
 Then,
with
probability $1-O(\frac{1}{q})$, we have the followings.

\mn
(a) There are at most $\frac{5q}{6}$  counterfeit coins $c$ that satisfy $|w(c)|< \ga$  (not exclusive) or belong to  a set $A_{0,j}$ containing more than one counterfeit coin, $j=1,...,2^{\lmq}$.

\mn
(b) For each $i=1,..., \lceil 2\log q \rceil-1$, $A_{i,j}$ contains
at most   $\frac{i+2\log q}{i}$ counterfeit coins.

\mn (c) For each  $i=1,..., \lceil 2\log q \rceil-1$, there are at most  $2^{-(i+1)}q +
q^{3/4}$ sets $A_{i,j}$ that
 contain more than one counterfeit coin.

\mn  (d) For $i\geq \lceil 2\log q\rceil-1$,   each $A_{i,j}$ contains
one or less counterfeit coin.

\mn (e) Each  $A_{\lceil 3\log n \rceil,j}$ contains at most one coin.
\end{lemma}

Now we are ready to present the algorithm described in Theorem \ref{cw}.

 \mn {\bf Algorithm }  (i) (Initially, $q=m$ and $A$ is the set of all $n$ coins.) Construct random subsets $A_{i,j}$ of $A$ as
above with parameter $q$. Then weigh $A_{0,j}$ for all $j=1,...,
2^{\ell_{q}}$, and denote $w_{0, j} = w(A_{0,j}) $ $j=1,...,2^{\ell_{q}}$ and $J_0$ to be the set of all $j$ such that
$|w_{0,j}|\geq \ga$. Then go to (ii), where, in general, $w(B)=\sum_{c\in B} w(c)$ for a set $B$ of coins.

\mn (ii) (Initially $ i=1$ and  $J=J_0$.) 
After relabeling, we may assume  $J=\{1, ..., |J|\}$.
Apply Corollary
\ref{sm} with $\gc_{_i} = \max\{\lceil\log (\frac{6\gb}{\ga}) \rceil, \lceil
\log (\frac{3\gb (i+2\log q)}{i\ga }) \rceil\}$ to $A_{i,2}, ..., A_{i,2|J|}$ and obtain $x_{r}$ satisfying
\begin{equation}\label{main3}
  w (A_{i,2r}) = x_{_{r}}  -
   \sum_{k=1}^{r-1} a_{_{rk}} w (A_{i,2k})  -
    \sum_{k=1}^{k_{r}} \frac{
  w(A_{i,2(r+k)}) }{2^{k\gc_{i}}}.
  \end{equation}

 \newcommand{\case}[4]{
\left\{ \begin{array}{ll} {#1} & \mbox{#2} \\ {#3}  & \mbox{#4}
\end{array} \right.}

\newcommand{\caseth}[6]{
\left\{ \begin{array}{lll} {#1} & \mbox{#2} \\ {#3}  & \mbox{#4}
\\ {#5}  & \mbox{#6}
\end{array} \right.}

\mn
Set, inductively in $r=1,..., |J|$,
 \begin{equation}\label{uu}
u_{_{2r}}=\caseth{w_{i-1,r}}{ if $|x_{_{r}} -\sum_{k=1}^{r-1}a_{
rk}u_{2k}|\geq \frac{\ga}{2}$}{}{}{~~~~0}{~otherwise,}
 \end{equation}
and $u_{_{2r-1}}= w_{i-1,r}-u_{_{2r}}$, $r=1,...,|J|$. Go to (iii) if $i< \lceil 2\log q \rceil$. Otherwise, go to (iv).

\mn
  (iii) (Initially, $s=-2$.)  Randomly select each $j$ satisfying $u_{\! _j}=0$ and  $j \leq \min\{s, 2|J|\}$ with probability $1/2$, independently of all other $j$. Weigh $\cup\{A_{i,j}:{\rm selected}~ j\} $. The weight is $0$ if no $j$ is selected. Do this random weighing $\lceil \log (i^2+1)\rceil+3$ times, independently of all other random weighings. This procedure is called a random test at $s$. If the test is passed, i.e., all weights are $0$, then  update $s$ to be $s+2i^2$.
  If it is failed and $s\leq 2|J|$, correct $u_{s}$ by weighing $A_{i, s}$, that is, weigh $A_{i, s}$, and  update  $u_{s}$  to be  $w(A_{i, s})$ and $u_{s-1}$ to be $w_{i-1, s/2}-u_s$. (Note that $s$ is even.) Update also $u_{j}$ for all  $j > s$ according to (\ref{uu}) and $u_{_{2r-1}}=w_{i-1,r}-u_{_{2r}}$.  If the test is failed and $s> 2|J|$, then do nothing.
  Update $s$ to be $s-2$ for both cases. This step including all updating
    is to be called a correction step of $u_{s}$, or simply a correction step, even for $s>|J|$. It does not necessarily mean that $u_{_s}$ was not $w(A_{i,s})$ just before the correction step though.

  \mn If $s\leq 2|J|+ 8i^2
  \log q$, repeat (iii) with updated $s$. Otherwise, let $w_{i,j}=u_j$, $j=1,...,2|J|$. Then return to the original label and update  $i$, $J$ to be $i+1$, $\{ j: \mbox{$w_{i,j}$ is defined}$ ${\rm and}~|w_{i,j}| \geq \ga \}$,   respectively, and go to (ii).

  \mn
  (iv) Set $w_{i,j}=u_j$, $j=1,...,2|J|$. Then return to the original label and update $J$ to be $\{ j: \mbox{$w_{i,j}$ is defined}$ ${\rm and}~ |w_{i,j}| \geq \ga \}$.   If $i< \lceil 3\log n \rceil $, then go to (ii) after updating $i$ to be $i+1$.
  If $i= \lceil 3\log n \rceil $, then output $J$ and declare that   all coins in $\cup_{j\in J} A_{i,j}$   are  counterfeit.
Remove all coins that are declared counterfeit from the set $A$ of  all coins and update $q$ to be $5q/6$. If $q> m^{0.8}+ 2\eps m $ go to (i). Otherwise, go to (v).

\mn
  (v) Apply the randomized binary
search to find counterfeit coins one by one,  using   $(\lceil \log n\rceil+3 )(m^{0.8}+2 \eps m)$ queries.

\mn

The core parts of the algorithm are (ii) and (iii).
If $w_{i-1,j}=w(A_{i-1,j})$ and every set $A_{i-1,j}$ contains
at most one coin, then $w(A_{i,2j})= 0$ or $w_{i-1,j}$.
Provided $|\sum_{k=1}^{k_{r}} \frac{
  w(A_{i,2(r+k)}) }{2^{k\gc_{i}}}|$ is small enough, say less than
  $\ga/2$ (see (a) lemma \ref{main}), it is not hard to show that $u_{2r}= w(A_{i,2r})$ and $u_{2r-1}= w(A_{i,2r-1})$ for all $r$. (See Corollary \ref{bi}.)
  This was one of main ideas of Bshouty and Mazzawi \cite{BM10_MFCS}.
  In general, as some sets $A_{i-1,j}$  contain more than one counterfeit coin, $u_{2r}$ may or may not be $w(A_{i,2r})$.

  If  $r$ is the smallest $r$ with $u_{2r}\not =w(A_{i,2r})$,
  $u_{2r'}=w(A_{i,2r'})$, $r' > r$, is not guaranteed any more
  even if the  set $A_{i-1,r'}$ contains only one counterfeit coin.
    This is why we introduced
  random tests and correction steps in (iii). The random tests
  generate a random walk that travels according to the value of $s$.
  It turns out that the walk goes forward until it passes or at $2r$.
  Once the random walk passes or at $2r$, it goes backward
  with a probability close  to $1$ (not extremely close to $1$ though).
  It is expected that the random walk with correction steps
  quickly identifies and corrects  $u_{_{2r}}$.

  Moreover, it turns out that $r$ is the smallest $r$
  with $u_{2r}\not =w(A_{i,2r})$ only if $A_{i-1,r}$ contains
   more than one counterfeit coin.
  If not many sets $A_{i-1,r}$ contain more than one counterfeit coin
  (see (c) of Lemma \ref{random}), the number of queries asked
  to identify and correct  corresponding $u_{2r}$'s seems
   to be reasonably small.
  In other words, the lesser the number
  of sets $A_{i-1,r}$ containing more than one counterfeit coin is,
  the faster $s$ increases. Eventually, $s$ keeps increasing after
  all corresponding $u_{2r}$'s
  are corrected.

\mn
{\bf Remark.} (a) Though the initial value $-2$ of $s$ looks somewhat strange,
it is natural as  $s=2r-2$ when $u_{2r}$ is corrected and
the initial value must be determined as if  the imaginary $u_0$ were corrected.

\mn
(b) When the random test fails, it may be tempting to  find $A_{i,2r}$
with $w(A_{i,2r})\not=0$, say using a binary search. However,
the number of queries needed to find such a set can be as large
as $\Omega(\log q)$, while our algorithm is expected
to correct $u_{2r}$ using $O(i^2 \log (i^2+1))$ queries.
This save queries when $i$ is small. Though the bound
is not extremely good if $i$ is large, it is not really a
matter as there are much less sets  $A_{i-1,j}$ containing
more than one counterfeit coin. (See (c) of Lemma \ref{random}.)

\mn

\mn
To analyze  the algorithm, we precisely summarize core properties of the algorithm.

\begin{lemma}\label{main} Suppose (a)-(e) of Lemma \ref{random} hold  for $q$ and $w_{i-1,j} = w(A_{i-1,j})$ for a fixed $i=1,..., \lceil 3\log n \rceil$ and all $j=1,..., |J|$.  Then  we have the followings.

\noindent
(a) For all $r=1,..., |J|$, $ \Big|\sum_{k=1}^{k_{r}} \frac{
  w(A_{i,2(r+k)}) }{2^{k\gc_{i}}}\Big| < \frac{\ga}{2}$.

\noindent
(b)  If $r$ is the smallest $r$ such that $u_{_{2r}} \not= w(A_{i,2r})$ when $u_{2r}$ is first defined or updated, then neither  $w(A_{i,2r-1})$ nor $w(A_{i,2r})$ is zero, especially $A_{i-1,r}$ contains more than one counterfeit coin.

\noindent
(c) Suppose  $i< \lceil 2\log q\rceil$ and   $u_{_j} = w(A_{i,j})$ for
all $j\leq 2r-2$ at a step.
If $s\leq 2r-2$ at the step, then $s$ keeps
increasing  until $s\geq 2r$. And once $s\geq 2r$, $s\geq  2r$
at all the following steps except possibly one step,
which is a correction step of $u_{_{2r}}$ and $s=2r-2$.

\old{\noindent
(d) Suppose $i< \lceil 2\log q\rceil$. If  $s\geq 2r$ and $u_{2r} \not= w(A_{i,2r})$
at  a step, then the probability that $s$ increases at the next step is at most $\frac{1}{8(i^2+1)}$.}
\end{lemma}

\proof
(a) For $i<\lceil 2\log q \rceil$, since $i \leq  2\log q $, $2^{^{\gc_i}} \geq \frac{3\gb(i+2 \log q)}{i\ga}$ and $|w(A_{i,2(r+k)})|\leq \frac{\gb(i+2\log q)}{i}$ (as $A_{i,2(r+k)}$ contains at most $\frac{i+2\log q}{i}$ counterfeit coins by (b) of Lemma \ref{random}),
  we have
$$ \Big|\sum_{k=1}^{k_{r}} \frac{
  w(A_{i,2(r+k)}) }{2^{k {{\gc_{i}}}}}\Big|
  \leq \sum_{k=1}^{k_{r}} \frac{
 \gb(i+2\log q) }{i(\frac{3\gb(i+2 \log q)}{i\ga})^k}
 \leq \frac{\ga}{3} +\frac{\ga}{3}  \sum_{k=1}^{\infty}
 \Big( \frac{2\ga \log q }{12\gb \log q} \Big)^k
 \leq \frac{\ga}{3} +
  \frac{\ga}{3} \sum_{k=1}^{\infty} \Big(\frac{1
  }{6 }\Big)^k < \frac{\ga}{2}.
$$
If $i\geq \lceil 2\log q \rceil$, then  each $A_{i,j}$ contains
one or less counterfeit coin by (d) of Lemma \ref{random}, which together with $2^{\gc_i} \geq \frac{6\gb }{\ga}$ gives
  $$ \Big|\sum_{k=1}^{k_{r}} \frac{
  w(A_{i,2(r+k)}) }{2^{k\gc_{i}}}\Big|
  \leq \sum_{k=1}^{k_{r}} \frac{
 \gb }{(\frac{6\gb }{\ga})^k}
 \leq \frac{\ga}{6} +\frac{\ga}{6}  \sum_{k=1}^{\infty}
 \Big( \frac{\ga}{6\gb } \Big)^k
 \leq \frac{\ga}{6} +
  \frac{\ga}{6} \sum_{k=1}^{\infty} \Big(\frac{1
  }{6 }\Big)^k < \frac{\ga}{2}.
$$

\mn
(b) As $r$ is the smallest $r$ such that $u_{_{2r}} \not= w(A_{i,2r})$ when $u_{_{2r}}$ is defined or updated, $u_{2j} = w(A_{i,2j})$ for all $j< r$ and hence
$$ w (A_{i,2r}) = x_{_{r}}  -
   \sum_{k=1}^{r-1} a_{_{rk}} w (A_{i,2k})  -
    \sum_{k=1}^{k_{r}} \frac{
  w(A_{i,2(r+k)}) }{2^{k\gc_{i}}}=  x_{_{r}}  -
   \sum_{k=1}^{r-1} a_{_{rk}} u_{_{2k}}  -
    \sum_{k=1}^{k_{r}} \frac{
  w(A_{i,2(r+k)}) }{2^{k\gc_{i}}}.$$
If $u_{_{2r}}=0$, $u_{_{2r}}\not=w (A_{i,2r})$ yields that $w (A_{i,2r})\not=0$. On the other hand, $u_{_{2r}}=0$ implies that
$|x_{_{r}}  -
   \sum_{k=1}^{r-1} a_{_{rk}} u_{_{2k}}| < \ga/2$. This together with (a) gives that $|w (A_{i,2r})| < \ga$. Since $|w (A_{i-1,r})|=|w_{i-1
   ,j}| \geq \ga$ and $w (A_{i,2r-1})= w (A_{i-1,r})-w (A_{i,2r})$, $w (A_{i,2r-1})\not=0$.
If $u_{_{\! 2r}}=w_{i-1,r}(=w(A_{i-1,r}))
$, then  $u_{_{\! 2r}}\not=w (A_{i,2r})$ yields that
$w (A_{i,2r})\not=w (A_{i-1,r})$ and hence $w (A_{i,2r-1})= w (A_{i-1,r})-w (A_{i,2r})\not=0$. On the other hand, $u_{_{\! 2r}}=
w_{i-1,r}$ implies that
$|x_{_{r}}  -
   \sum_{k=1}^{r-1} a_{_{rk}} u_{_{2k}}| \geq \ga/2$. This together with (a) gives that $|w (A_{i,2r})| >0$, i.e., $w (A_{i,2r})\not=0$

\mn (c) We prove this by reverse induction. For $r=|J|+1$,
if $u_{j}= w(A_{i,j})$ for all $j\leq 2r-2=2|J|$,  then $w(A_{i,j})=0$ whenever $u_{j}=0$, for all $j\leq 2|J|$. Thus, the random test must be passed and   $s$ keeps increasing regardless of the value of $s$ (as no $u_{_j}$ is updated).
Suppose $u_{j}= w(A_{i,j})$ for all $j\leq 2r-2$ with $r\leq |J|$. Then  $w(A_{i,j})=0$  for all $j\leq 2r-2$ with $u_{j}=0$. If $s \leq 2r-2$, the random test  must be passed and   hence $s$ increases. Once $s>2r-2$, or equivalently $s\geq 2r$ (as $s$ is even), no $u_{j}$ with $j\leq 2r-2$ is updated before a correction step of $u_{_{\! 2r}}$.
If there is no correction step of $u_{_{\! 2r}}$, then $s\geq 2r$ at all the following steps. If $u_{_{\! 2r}}$ is corrected at a step,  then $s=2r-2$ and $u_{j}= w(A_{i,j})$ for all $j\leq 2r$ at the step. The induction hypothesis especially yields $s\geq 2r$ at all  steps after the correction step.

\mn
\qed

\renewcommand{\S}{{\cal S}}
\mn

To analyze (iii) of the algorithm for a fixed
$i<\lceil 2\log q\rceil$, we may regard the whole process
as a random walk $\S$ that travels
according to the value of $s$. That is, $\S=(s_{_0}, s_{_1}, ...)$,
where $s_{_{k}} $ is the value of $s$ at the end of the $k^{\rm th}$ step.
Note that $\S$ goes backward, i.e.,  $s$ decreases,
at a step if and only if the step is a correction step.
We will see that $\S$ goes forward until it passes or at $2r$ for the
the smallest $r$ with  $u_{_{2r}}\not=w(A_{i,2r})$, and
then $\S$ tends to go backward
with probabilities  close (not extremely though) to $1$ until
$u_{_{2r}}$  is corrected.

We  partition $S$ into  a few subrandom walks that are essentially
independent identically distributed (i.i.d).
They are not exactly i.i.d though.
Let $r_{_{\! 0}}=0$. The $0^{\rm th}$ (sub)random walk $\S_0$ (of $\S$)
starts when the whole process starts  and ends at the same time,
that is, $\S_0= (s_{_0})$ (recall $s_{_0}=-2$).
Let $r_{_{\! 1}}$ be the
the smallest $r$ with  $u_{_{2r}}\not=w(A_{i,2r})$.
The first random walk $\S_1$
starts immediately after  $\S_0$ ends and it ends when
$s=2r_{_{1}}-2$ at a backward step  or  $s> 2|J|+ 8i^2
  \log q$ for the first time, whichever comes first.
  Generally, for $\ell\geq 2$, if $\S_{\ell-1}$ ends with  $s=2r_{_{\el-1}}-2$,
  then let $r_{_{\el}}$ be the smallest $r\leq |J|$
  such that $u_{_{2r}} \not= w(A_{i, 2r})$ at the step $\S_{\ell-1}$ ends.
   The $\ell^{\rm th}$ random walk $\S_\ell$ starts immediately after $\S_{\ell-1}$
   ends, and it ends when
$s=2r_{_{\el}}-2$ at a backward step or  $s> 2|J|+ 8i^2
  \log q$ for the firs time, whichever comes first.
  However, $\S_{\ell}$ does not end at a
  forward step with $s=2r_{_{\el}}-2$. In theory,
  it is possible that $\S_{\ell}$ is infinite, though
 it is not difficult to show that the probability of $\S_{\ell}$
  being infinite is $0$.
Both of $r_{_{\el\,'}}$ and $\S_{\ell\,'}$ are not defined for all $\ell' \geq \ell$,
if $\S_{\ell-1}$ is infinite or ends with $s> 2|J|+ 8i^2
  \log q$, or $u_{_{2r}} = w(A_{i, 2r})$ for all $r\leq |J|$ at the last step of $\S_{\ell-1}$.

  The random walk $\S_{\ell}$ is called good if it is defined
  and ends with $s=2r_{_{\el}}-2$.
     Note that the last step of a good random walk $\S_\ell$  is
  the first correction step of $u_{_{2r_{_{\el}}}}$ after $\S_\ell$ starts.
  In other words, a good random walk $\S_{\ell}$ ends when it corrects
  $u_{_{2r_{\el}}}$ where $r_{\el}$ is the smallest $r$ such
  that $u_{_{2r}} \not= w(A_{i, 2r})$ when it starts.
  We also note that  $r_{_{\el}}$, $\S_{\ell}$ are defined only
  if $\S_{\el-1}$ is good. In particular,  $\S_{\ell}$ is good only if
 $ \S_{\ell\, '}$ is good for all $\ell' \leq \ell-1$.

\old{
last step of $\S_{\ell-1}$  If there is no correction step of $u_{_{2r_{_{\el-1}}}}$ after $r_{_{\el-1}}$ is defined, , then $r_{_{\el}}$ is not defined. The $\ell^{\rm th}$ random walk $\S_\ell$ starts from (the end of) the step at which $s\geq 2r_{_\ell}$ for the first time after $r_{_{\el}}$ is defined and it ends at the step at which $s=2r_{_{\el}}-2$ or $s> 2|J|+ 8i^2
  \log q$ for the first time after it starts. Note that, if the (sub)random walk $\S_\ell$ ends with $s=2r_{_{\el}}-2$, then $u_{_{2r_{_\el}}}$ is corrected at the last step of $\S_\ell$, and then $r_{_{\el+1}}$ is defined at the last step.
  \old{By convention, the $0^{\rm th}$ step means when (iii) starts. The $0^{\rm th}$ random walk $\S_0$ starts with $s=0$ at  the $0^{\rm th}$ step with $s=0$  and ends}

The above lemma basically says that }

\mn
\begin{corollary}\label{m2} Under the same hypotheses as in Lemma \ref{main} with $i\leq \lceil 2\log q\rceil -1$, we have the followings.

\mn (a) For $\ell\geq 1$,  if  $\S_{\ell}$ is good, $u_{j}= w(A_{i,j})$
for all $1\leq j\leq 2r_{_{\el}}$ at the last step of $\S_\ell$,
especially $r_{_{\el+1}}> r_{_{\el}}$ if $r_{_{\el+1}}$ is defined.
Furthermore, a good random walk $\S_{\ell}$ starts
with $s=2r_{_{\el-1}}-2$ and keeps going forward
until $s\geq 2r_{_{\el}}$, and then goes back and force
with $s\geq 2r_{_{\el}}$ at all steps except the last step
at which $s= 2r_{_{\el}}-2$.

\mn
(b) If $r_{_\el}$ is defined, then neither  $w(A_{i,2r_{_\el}-1})$ nor $w(A_{i,2r_{_\el}})$ is zero and  $A_{i-1,r}$ contains more than one counterfeit coin.  In particular, $r_{_{\el}}$ and $\S_\ell$ are not defined if $\ell > h_q:=\lceil 2^{-(i+1)}q+q^{3/4}\rfloor$.

\mn
(c) Suppose every $\S_{\ell}$  is good if defined. Then $w_{i,j}= w(A_{i,j})$, for all $j=1,...,2|J|$.
\old{and, for $r=1,..., |J|$, $r=r_{_{\el}}$ for some $\ell\geq 1$ if and only if neither  $w(A_{i,2r-1})$ nor $w(A_{i,2r})$ is zero.}

\end{corollary}

 \proof
(a)   For $\ell\geq 1$, suppose $u_{j}= w(A_{i,j})$ for all
$j\leq 2r_{_{\el-1}}$ at  the last step of  $\S_{\ell-1}$.
(This is trivial when $\ell=1$.) Since $\S_{\ell}$ is good only
if $\S_{\ell-1}$ is good,  the induction hypothesis may be applied to
obtain   $r_{_{\el}}> r_{_{\el-1}}$ and hence
$$u_{j}= w(A_{i,j})~~\mbox{ for all}~~  j\leq 2r_{_{\el}}-2$$
 at the last step
of $\S_{\ell-1}$. Then, (c) of Lemma \ref{main}  gives that $s$ keeps
increasing (without updating $u_j$) after  the last step of $\S_{\ell-1}$,
at which $s= 2r_{_{\el-1}}-2$, until $s\geq  2r_{_{\el}}$.
Once $s\geq  2r_{_{\el}}$, no $u_j$ with $j\leq 2r_{_{\el}}-2$ is updated
before the last step of $\S_\ell$.
Since $\S_\ell$ is good, $u_{2r}$ is corrected and hence
 $u_{2r}= w(A_{i,2r})$, $u_{2r-1}= w(A_{i,2r-1})$
at the last step of $\S_\ell$. The second part is already shown too.

\mn
(b) Since  $r_{_{\el}}> r_{_{\el-1}}$ is defined, $s= r_{_{\el-1}}-2$ at
the last step  of $\S_{\ell-1}$
and $u_{_{2r_{_{\el}}}}$ is updated at the last step of
$\S_{\ell-1}$. By (b) of Lemma \ref{main}, neither  $w(A_{i,2r-1})$
nor $w(A_{i,2r})$ is zero.
The second part follows from that
 all $r_{_{\el}}$'s are distinct and there are at most
 $h_q$ sets $A_{i-1,r}$ containing more than one
 counterfeit coin (see (c) of Lemma \ref{random}).

\mn
(c) For  the largest $\ell$ for which $r_{_\el}$ is defined, as  $\S_\ell$ is good and $r_{_{\el+1}}$ is not defined,
 $u_{_{2r}}=w(A_{i,2r})$  for all $r=1,..., |J|$ at the last step of $\S_\ell$.
Since  $s$ keeps
increasing after the last step  and eventually $s> 2|J| + 8i^2\log q$ without updating $u_{_{j}}$'s, we have $w_{i,2r} = w(A_{i,2r})$ for $r=1,...,|J|$, and $w_{{i,2r-1}}=w_{i-1,r}- w_{i,2r}=w(A_{i-1,r})-w(A_{i,2r})=w(A_{i,2r-1})$.

\old{For the second part, it is enough to show that $r=r_{\ell}$ for some $\ell$ provided that neither  $w(A_{i,2r-1})$ nor $w(A_{i,2r})$ is zero, as  the other direction is proven in (b).
We prove the contrapositive of the statement. Suppose $r\not=r_{_{\el}}$ for all $\ell\geq 1$ and $\ell$ is   the largest $\ell \geq 0$ such that  $r_{_{\el}} \leq r$. Then
 $r> r_{_{\el}} $ and   $u_{_{2r}}$ is updated at the last step of $S_{\ell}$, for $\S_\ell$ is good. Thus, $u_{_{2r}}$ is either $0$ or $w_{i-1,r}=w(A_{i-1,r})$ at the last step.
  If $r_{_{\el+1}}$ is defined, then $r_{_{\el+1}}>r$ and $u_{_{j}}= w(A_{i,j})$ for all $j\leq 2r_{_{\el+1}}-2$, especially $u_{_{2r}}= w(A_{i,2r})$ at the last step of $S_{\ell}$. If $r_{_{\el+1}}$ is not defined, $u_{_{j}}= w(A_{i,j})$ for all $j\leq 2|J|$, especially $u_{_{2r}}= w(A_{i,2r})$ at the last step of $S_{\ell}$.
  However, as $u_{_{2r}}$ is either $0$ or $w(A_{i-1,r})$,
  we have $w(A_{i,2r})=w(A_{i-1,r})-w(A_{i,2r-1})=0$ or $w(A_{i-1,r})$,
 which implies $w(A_{i,2r})=0$ or $w(A_{i,2r-1})=0$. }

\mn \qed

\mn

If $\lceil 2 \log q \rceil\leq i\leq  \lceil 3\log n \rceil$,
each $A_{i,j}$ contains  at most one counterfeit coin by (d) of Lemma \ref{random}. Then it easily follows that $w_{i,j}=w(A_{i,j})$ for all $j=1,..., 2|J|$.

\begin{corollary}\label{bi} Under the same hypotheses as in Lemma \ref{main}
with $\lceil 2 \log q \rceil\leq i\leq  \lceil 3\log n \rceil$,
$w_{i,j}=w(A_{i,j})$ for all $j=1,..., 2|J|$.
\end{corollary}

\proof Recall  that $w_{i,j}=u_{_j}$, where $u_{_j}$'s
are defined in (ii). Take, if any, the smallest $r$ such
that $u_{_{2r}}\not=w(A_{i,2r})$. Then, (b) of Lemma \ref{main}
implies that $A_{i-1,r}$ contains more than one counterfeit coin,
which is not possible as each $A_{i-1,r}$ contains  at most one counterfeit
coin due to (d) of Lemma \ref{random}. Hence, $w_{i,2r}=w(A_{i,2r})$
for all $r=1,..., |J|$ and
$w_{i,2r-1} = u_{_{2r-1}}=w(A_{i-1,r})-w(A_{i,2r})=w(A_{i,2r-1})$.

\mn \qed

Corollaries \ref{m2} and \ref{bi} provide all but one basic properties to analyze  the algorithm.
The missing property is that, with high probability, every $\S_\ell$ is good if defined,  the hypothesis of (c) of Corollary \ref{m2}.
For  the query complexity,  an upper bound for the number $|\S_\ell|$ of steps in  $\S_\ell$ is needed. As we want to bound $|\S_\ell|$ only for good $\S_{\ell}$, we will consider  $|\S_\ell|\chi_{_\ell} $, where
$$ \chi_{_\ell} =\caseth{1}{if $\S_\ell$ is good}{}{}{0}{otherwise.}$$

It will be first shown that, after $\S_\ell$ passes or at $2r_{_\el}$,
the random walk goes to backward with probability at least
$1-\frac{1}{8(i^2+1)}$ until $u_{_{2r_{_\el}}}$ is corrected,
which follows from  $w(\cup\{A_{i,j}:{\rm selected}~ j\}) \not=0$
with probability at least $1/2$ during the process.
Thus, $\S_\ell$
goes  backward by at least $7/4$ in expectation
after $\S_\ell$ passes or at $2r_{_\el}$,
as $\S_\ell$ goes forward by $2i^2$ with probability at most
$\frac{1}{8(i^2+1)}$ and goes backward by $2$ otherwise.
This is why  $\S_\ell$ is expected to be  good.
The number $F_\ell$ of forward steps in $\S_\ell$ after it passes
or at $2r_{_\el}$ is also expected to be reasonably small,
namely $O(1)$ with a probability close to $1$.
It actually turns out that the probability of $F_\ell=k$ is
at most $e^{-k+1}$ and the sum $\sum_{\ell=1}^{h_q}  F_\ell$
may be bounded by $O(h_q)$ with high enough probability, say
with probability $1-e^{-\Omega(q^{3/4})}$, where
$h_q=\lceil 2^{-(i+1)}q+q^{3/4}\rfloor$ as in (b) of Lemma \ref{m2}.

Then, it is not difficult to show that
the number of all steps in $S_{\ell}$ after it passes or at $2r_{_\el}$ is $O( i^2F_\ell)$, especially there are $O(i^2F_\ell)$ backward steps in $S_{\ell}$ by (a) of Corollary \ref{m2}. Therefore, there are $O(i^2 h_q)$ backward steps in $\S$ with high probability. All other steps in $\S$ are forward steps and hence there are
$$O\Big(i^2 h_q+ \frac{|J|+8i^2 \log q+2i^2h_q}{2i^2} \Big)
=O\Big (\frac{|J|}{i^2}+ (i^2+2)\Big(\frac{q}{2^{i+1}}+q^{3/4}\Big)\Big)$$
steps in $\S$. As $O(\log (i^2+1))$ queries are asked at each step, the number of queries asked in the $i^{\rm th}$ round, $i=1,...,\lceil 2\log q \rceil-1$, is $ O\Big(q\Big(\frac{1}{i^2}+ \frac{i^2 }{2^i}\Big)\log (i^2+1) \Big)$ assuming $|J|\leq q$.

The precise statements are presented in the next lemma. Though  idea is simple as illustrated above, our proof of the lemma is somewhat lengthy, partly because it is proven  rigorously without referring other theories. We prove the lemma at the end of this section. Readers familiar with random walks may skip the proof.

\begin{lemma}\label{m5s} Under the same hypotheses as in Lemma \ref{main} with $i\leq \lceil 2\log q\rceil -1$, if  $u_{2r} \not= w(A_{i,2r})$ and  $s\geq 2r$ at  a step, then the probability that $s$ increases at the next step is at most $\frac{1}{8(i^2+1)}$.
Moreover, every $ \S_\ell$ is good if defined, with probability $1-O(q^{-3})$, and the number $|\S|$ of all steps satisfies
$$  \pr \Big[ |\S| \geq \frac{|J|}{i^2}+ 4(i^2+2)\Big(\frac{q}{2^{i+1}}+q^{3/4}\Big)\Big]
 = O( q^{-3}).
$$

\end{lemma}

\old{
$$ \int_{0}^{\infty} \frac{(x^2+2)\log (x^2+1
)}{2^x} dx = \Big(-\frac{(x^2+2)\ln(x^2+2)}{2^x \ln 2} - \frac{2x\ln(x^2+2)}{2^x (\ln 2)^2} - \frac{2}{2^x (\ln 2)^3}\Big)\Big|_{0}^{\infty}. $$}

\mn

\mn
{\bf Correctness of the algorithm} Once Lemmas \ref{random}, \ref{m5s} and Corollaries \ref{m2}, \ref{bi} are established, it is easy to see that the algorithm finds counterfeit coins as desired. In the next lemma, we precisely describe it along with a property needed to bound query complexity.

\mn

\newcommand{\sub}{\subseteq}
 \renewcommand{\sup}{\supseteq}
\begin{lemma}\label{m4} For a fixed $q> m^{0.8}+ 2\eps m$, the followings hold
with probability
 $1-O(1/q)$, assuming the same in the prior round.

 \mn
 (a)  The statements (a)-(e) of Lemma \ref{random} hold.

 \mn
 (b) Whenever $w_{i,j}$ is defined, $w_{i,j}= w(A_{i,j})$. In particular,  a  coin   declared to be counterfeit must be counterfeit.

\old{
 \mn
 (c) The algorithm  defines  $w_{i,j}$  for every pair $i,j$ such that $A_{i,j}$ contains only one counterfeit coin $c$, and  $|w(c)|>\ga$ and $ A_{i,j}\sub A_{0,\ell}$ for  $A_{0,\ell}$ containing no counterfeit coin other than $c$.}

\mn
 (c) The algorithm finds  every counterfeit coin $c$ with $|w(c)|\geq \ga$ that is a unique counterfeit coin  of  $A_{0,\ell}$ for some $\ell=1,...,2^{\q}$. And the number of remaining counterfeit coins is at most  the updated $q$.

 \mn
(d) The number of queries asked in all rounds of (iii) is
$O(q)$, where the constant in $O(q)$ is at most
$\sum_{i=1}^{\infty} \frac{5+\log (i^2+1)}{i^2} +\frac{(i^2+2)(5+\log (i^2+1))}{2^{i-1}}+o(1).$

 \end{lemma}

\mn
\proof As $q> 2\eps m$ and (c) holds in the prior round, Lemma \ref{random} yields that the statements in (a) hold with probability $1-O(1/q)$. We assume the statements to prove the other properties.

To prove  the other properties, we further assume that every $\S_\ell$ is good if defined and that, for each $i=1,..., \lceil 2\log q \rceil-1$ and  the number $|\S|$ of all steps in the $i^{\rm th}$ round of (iii),
 \beq\label{ss} |\S| \leq \frac{|J|}{i^2}+ 4(i^2+2)\Big(\frac{q}{2^{i+1}}+q^{3/4}\Big),
  \enq
  both of which hold with probability $1-O(q^{-3})$ by Lemma \ref{m5s}. Then the first part of (b) follows from (c) of Corollary \ref{m2}, and Corollary \ref{bi}. Since every $A_{\lceil 3\log n \rceil,j}$ consists of one or no coin,
each coin $c$ in $\cup_{j\in J} A_{\lceil 3\log n \rceil,j}$
satisfies $|w(c)|= |w(A_{\lceil 3\log n \rceil,j})|=
|w_{\lceil 3\log n \rceil,j}|\geq \ga $ for some $j\in J$, especially, $c$ is counterfeit.

If a coin $c$ with $|w(c)|\geq \ga$ is a unique counterfeit coin in
$A_{0,\ell_0}$, then, for each $i=0,..., \lceil 3\log n \rceil$,
there is a unique $\ell_i$ such that  $A_{i, \ell_i}\sub A_{0,\ell_0}$
contains   $c$.  It is clear, by the way how  $A_{i,j}$'s are constructed,  that
$A_{i,\ell_i} \sub A_{i-1,\ell_{i-1}}$ for all $i=1,..., \lceil 3\log n \rceil$. Moreover, since $c$ is  a unique counterfeit coin  of  $A_{i,\ell_i}$, $|w(A_{i,\ell_i})|\geq \ga$ for all $i=0,..., \lceil 3\log n \rceil$. For $i=0$, $\ell_0\in J$ as $|w_{0,\ell_0}| = |w(A_{0,\ell_0})|\geq \ga$. For $i\geq 1$, assuming $\ell_{i-1} \in J$ when the prior round ends, $w_{i,\ell_i}= w(A_{i,\ell_i})$ by (c) of Corollary \ref{m2}, as $\ell_{i-1} \in J$ and $A_{i,\ell_i} \sub A_{i-1,\ell_{i-1}}$. Thus, $|w_{i,\ell_i}|=| w(A_{i,\ell_i})|\geq \ga$ implies that  $\ell_i$ is in the updated $J$.
We have just shown that $\ell_{i} \in J$ when the $i^{\rm th}$ round ends for each $i$, especially, for $i= \lceil 3\log n \rceil$. As $c
\in A_{i, \ell_i}$ for $i= \lceil 3\log n \rceil$, $c$ is declared
to be counterfeit. The second part of (c) follows from (a) of
Lemma \ref{random}.

Note that $|J|=|\{ j: \mbox{$w_{i,j}$ is defined}$
${\rm and}~|w_{i,j}| \geq \ga \}| \leq q$ by the  second part of (c) in the prior round and first part of (b), as $|w_{i,j}|=|w(A_{i,j})|\geq \ga$ implies that $A_{i,j}$ contains a counterfeit coin and the number of such sets is at most the number of counterfeit coins. Since  the algorithm asks at most $5+\log (i^2+1)$ queries at each step of $\S$ (one more query is needed in backward steps), \raf{ss} yields that
the number of queries  is at most
 $$ \sum_{i=1}^{\lceil 2\log q\rceil-1}\Big( \frac{(5+\log (i^2+1))q}{i^2}+ 4(i^2+2)(5+\log (i^2+1))\Big(\frac{q}{2^{i+1}}+q^{3/4}\Big)\Big)= O(q).
 $$

\mn
\qed


\mn

\mn

\mn

The lemma especially says that the number of remaining counterfeit coins decreases by factor $5/6$, with probability $1-O(1/{q})$. Applying this inductively until $q\leq  m^{0.8}+ 2\eps m$, we know the algorithm find all but at most $m^{0.8}+2\eps m$ counterfeit coins before it goes to (v), with probability $1-O(1/{m^{0.8}})$. All the remaining $m^{0.8}+2\eps m$ counterfeit coins
are found in (v), with probability $1-e^{-\Omega(m^{0.8})}$,  by Lemma \ref{rbs}.

\begin{corollary}\label{first}  The algorithm find all but at most $m^{0.8}+2\eps m$ counterfeit coins before it goes to (v), with probability $1-O(1/{m^{0.8}})$. All the remaining $m^{0.8}+2\eps m$ counterfeit coins
are found in (v), with probability $1-e^{-\Omega(m^{0.8})}$,  by Lemma \ref{rbs}.
\end{corollary}

 \mn
 {\bf Query Complexity} Suppose (a)-(d) of Lemma \ref{m4} hold for all $q$, which occurs   with probability $1-O(1/{m^{0.8}})$. Then
for each $q$,  the number of  remaining counterfeit coins is at most $q$. Especially, $|J|\leq q$ as seen in last paragraph of  the proof of Lemma \ref{m4}.
For each $q$,   $2^{\q} \leq 2q$ queries are needed in (i).
 For each $q$ and $i$, the number of queries asked in (ii) is
 $$ \frac{(2+o(1))\gc_{_i} |J|}{\log (\gc_{_i} |J|)}\leq  \caseth{
\frac{ (2+o(1))|J|}{ \log |J|}\Big\lceil \log (\frac{3\gb(i+2 \log q)}{i\ga})\Big\rceil}{
if $i< \lceil  2\log q\rceil$}{}{}{ \frac{ (2+o(1))|J|}{ \log |J|}
\Big\lceil\log (6\gb/\ga) \Big\rceil }{~if $i\geq  \lceil  2\log q\rceil$.}
 $$
Since $|J|\leq q$ and $$ \sum_{i=1}^{\lceil 2\log
q \rceil-1} \Big\lceil \log (\frac{3\gb(i+2 \log q)}{i\ga}) \Big\rceil
\leq  4\log
q  \log (3\gb/\ga)+ \log {2\lceil 2\log q \rceil-1 \choose \lceil 2\log q \rceil-1} \leq 4\log
q \log (3\gb/\ga)+4\log q + 1,
$$
and
$$\sum_{i=\lceil 2\log q \rceil}^{\lceil 3\log n \rceil }
\lceil\log (6\gb/\ga) \rceil  \leq 3\log (6\gb/\ga) \log n+ 3\log n
$$
for each $q$, the number of queries asked in (ii) is
$
 O\Big(\frac{q\log (\gb/\ga) \log n}{\log q}\Big)
. $

The number of all queries asked in (iii) for each $q$ is $O(q)$ by (d) of Lemma \ref{m4}.
No query is asked in (iv) and hence the total number of queries asked for fixed $q> m^{0.8}+2\eps m$ is $O\Big(\frac{q\log (\gb/\ga) \log n}{\log q}\Big)$.
As $q$ keeps decreasing  by factor of $5/6$, the number of queries asked before the algorithm goes to (v) is $O\Big(\frac{m\log (\gb/\ga) \log n}{\log m}\Big)$.

 \mn
 \qed

\bn

This together with Corollary \ref{first} implies that, if we artificially stop the algorithm when it asks $\frac{\eta m\log (\gb/\ga) \log n}{\log m}$ queries,  for the constant $\eta$ in the $O\Big(\frac{m\log (\gb/\ga) \log n}{\log m}\Big)$ term,
    all but at most $m^{0.8} +2\eps m$  counterfeit coins are found with probability $1-O(1/{m^{0.8}})$. As $(\lceil \log n\rceil+3 )(m^{0.8}+2 \eps m)$ queries are asked in (v) of the algorithm,
 Theorem \ref{cw} follows. We conclude this section by proving Lemma \ref{m5s}.

\bn
{\bf Proof of Lemma \ref{m5s}} \,  Note that each $u_{2r}$ may have one of
three values
 $0, w_{i-1, r }$, $w(A_{i,2r})$. Since $u_{2r} \not= w(A_{i,2r})$,
 $u_{2r}$ is either $0$ or $w_{i-1,j}$. If $u_{2r}=0$,
 then $w(A_{i,2r})\not=0$. If $u_{2r}=w_{i-1,j}(=w(A_{i-1,j}))$,
 then $u_{2r-1}=0$ while $w(A_{i,2r})\not= u_{2r}=w(A_{i-1,j})$
 yields $w(A_{i,2r-1})=w(A_{i-1,j})- w(A_{i,2r})\not=0$.
Particularly, there is $\ell\leq s$ such that $w(A_{i,\ell})\not=0$
while $u_{\ell} =0$. Suppose the random selection other than $\ell$ is
carried out. Then
the set of coins to be weighed is either
$\cup\{A_{i,j}:{\rm selected}~ j, j\not= \ell \}$ or
$\cup\{A_{i,j}:{\rm selected}~ j, j\not= \ell \} \cup A_{i,\ell}$,
each of which occurs with probability $1/2$.
Since $w(A_{i,\ell})\not=0$ implies that the weights of the two sets are
different,
$$ \pr [w( \cup \{A_{i,j}:{\rm selected}~ j \})=0] \leq 1/2. $$
After independently performing this $\lceil \log (i^2+1)\rceil+3$ times,
the probability that all weights are $0$ is at most
$2^{\lceil \log (i^2+1)\rceil+3}\leq  \frac{1}{8(i^2+1)}$.
That is, $s$ increases at the next step with probability at
most $\frac{1}{8(i^2+1)}$.

For the second part, suppose $\S_\ell$  is defined but it is not good,
which especially means that $\S_{\ell-1}$
is good. Then $\S_\ell$  must be infinite or reach a step with
$s> 2|J|+8i^2 \log q$.
As $\S_{\ell}$ starts with $s=2r_{_{\el-1}}-2$, $r_{_{\el-1}} < r_{_{\el}}$, and $u_{2r}= w(A_{i,2r})$ for all $r\leq 2r_{_\el}-2$,  the random walk $\S_\ell$   keeps going forward until $s\geq 2r_{_{\el}}$ by (c) of Lemma \ref{main}. Let  $\gs_{_\el}$  be the value of $s$ when  $\S_\ell$ reaches a step with $s\geq 2r_{_{\el}}$ for the first time. Then
  \beq\label{tl}
             2r_{_{\el}}\leq \gs_{_\el} \leq 2r_{_{\el}}+2i^2-2,~~{\rm or} ~~0\leq \gs_{_\el}/2 -r_{_{\el}}\leq i^2-1,
  \enq
  for $s$ increases by $2i^2$.
Hence,  there must be at least $\lfloor 4\log q\rfloor$ more forwarding steps to reach a step with $s> 2|J|+8i^2 \log q$, as,  otherwise,
$$s\leq \gs_{_\el}+2i^2 (\lfloor 4\log q\rfloor-1)
\leq 2r_{_{\el}}-2+2i^2 +2i^2 (\lfloor 4\log q\rfloor-1)\leq 2|J|+8i^2 \log q. $$
If $\S_\ell$ is infinite, there must be at least $\lfloor 4\log q\rfloor$ more forwarding steps too.

Counting after $\S_\ell$ reaches a step with $s\geq 2r_{_{\el}}$ for the first time, let $T$ be the number of steps in $\S_\ell$ until there are $\lfloor 4\log q\rfloor$ more forwarding steps.
For $\S_\ell$ is not good, there is no correction step of $u_{_{2r_{_\el}}}$, or equivalently $s \geq 2r_{_\el}$ after the count starts, particularly,   $T$  satisfies
$$  \gs_{_\el} + 2i^2 \lfloor 4\log q\rfloor -2(T-\lfloor 4\log q\rfloor )\geq 2r_{_{\el}},$$
which, together with \raf{tl}, gives
$$
 T \leq (i^2+1)\lfloor 4\log q\rfloor +  \gs_{_\el}/2 -r_{_{\el}}
 \leq (i^2+1)(\lfloor 4\log q\rfloor +1)
 . $$
 We have just shown that, for $t= (i^2+1)(\lfloor 4\log q\rfloor +1)$, \beq\label{ng}\pr [ \S_{\ell}~\mbox{is not good} ]
 \leq \pr \Big[~\exists\,  \lfloor 4\log q\rfloor
 ~\mbox{forward steps among the
 first $t$ or less steps of $\S_\ell$}\, \Big].  \enq
 To bound the last probability, it is convenient to introduce an auxiliary random walk $\S^*_\ell$.  The infinite random walk $\S^*_\ell$ starts when $\S_\ell$ reaches a step with $s\geq  2r_{_{\el}}$ for the first time and it is the same as $\S_\ell$ until $\S_\ell$ ends. Once $\S_\ell$ ends,  $\S^*_\ell$ keeps going forward by $2i^2$ with probability $\frac{1}{8(i^2+1)}$ and  backward by $2$
with probability $1-\frac{1}{8(i^2+1)}$.
Then, at any step, $\S^*_\ell$ goes forward with probability at most
$\frac{1}{8(i^2+1)}$.
\old{: If $S_{\ell}$ did not end before the step, then
when the step starts, $s\geq 2r_{\el}$ while $u_{_{2r_{_\el}}}\not= w(A_{i, 2r_{_\el}})$ remains unchanged. By the first part,  $\S_\ell$ goes forward at the step with probability at most $\frac{1}{8(i^2+1)}$ and so does $\S^*_\ell$.
If $S_{\ell}$ ended before the   step starts, the statement
clearly holds by the definition of $\S^*_\ell$.}

As  there are $\lfloor 4\log q\rfloor$ forward steps among the
 first $t$ steps of $\S^*_\ell$ if there are $\lfloor 4\log q\rfloor$ forward steps among the
 first $t$ or less steps of $\S_\ell$, \raf{ng} gives
 $$   \pr [ \S_{\ell}~\mbox{is not good} ]
 \leq  \pr \Big[ ~\exists\,  \lfloor 4\log q\rfloor
 ~\mbox{forward steps among the
 first $t$ steps of $\S^*_\ell$}\, \Big],$$
  which is at most ${t \choose \lfloor 4\log q\rfloor }
  \Big( \frac{1}{8(i^2+1)} \Big)^{\lfloor 4\log q\rfloor}$.
 Therefore,  using ${t \choose k} \leq (\frac{et}{k})^k$,
$$
\pr [ \S_{\ell}~\mbox{is not good} ]
\leq {t \choose \lfloor 4\log q\rfloor }
  \Big( \frac{1}{8(i^2+1)} \Big)^{\lfloor 4\log q\rfloor}
 \leq \exp\Big(  \lfloor 4\log q\rfloor \ln  \frac{e(i^2+1)(\lfloor 4\log q\rfloor+1)}{8(i^2+1)\lfloor 4\log q\rfloor}\Big).  $$
Using $\ln(e/8)\leq -1$ and $\ln (1+y) \leq y$ for $y\geq 0$,
we obtain
$$
\pr [ \S_{\ell}~\mbox{is not good} ]\leq
\exp\Big(-\lfloor 4\log q\rfloor+1\Big)=O(q^{-4}) .$$
Since $\S_\ell$ is defined for at most $h_q$ indices $\ell$ by (b) of Corollary \ref{m2},  and $h_q= O(q)$,  Boole's inequality yields the desired bound.

For the last bound, if $S_\ell$ is good,   let $F_\ell$ be  the number of all forward steps in $\S_\ell$ after $\S_\ell$ reaches a step with $s\geq 2r_{_{\el}}$ for the first time. If $S_\ell$ is not good or not defined, then $F_\ell=0$.
If $F_\ell=k \geq 1$, then $\S_\ell$ must be good and, for
the number $t$ of  all steps  in $\S_\ell$ after $\S_\ell$ reaches a step with $s\geq 2r_{_{\el}}$ for the first time, we have
$$ \gs_{_\el} + 2i^2 k -2(t -k)= 2r_{_{\el}}-2~~{\rm or} ~~t= (i^2+1)k + \gs_{_\el}/2-r_{_\el}+1\leq
 (i^2+1)(k+1), $$
(recall that $\gs_{_\el}$ is the value of $s$ when $\S_\ell$ reaches a step with $s\geq 2r_{_{\el}}$ for the first time).

 \old{and $F_\ell$ be  the number of forward steps in $\S_\ell$, both  until $\S_{\ell}$ ends. If $\S_\ell $ is not defined, $\gs_{_\el}=F_\ell=0$ Then, as $\S_\ell $ ends with $s=2r_{_{\el}}-2$,  $\gs_{_\el}$ and $F_\ell$ satisfy
 $$ \gs_{_\el} + 2i^2 k -2(t -k)= 2r_{_{\el}}-2~~{\rm or} ~~t= (i^2+1)k + \gs_{_\el}/2-r_{_\el}+1\leq
 (i^2+1)(k+1), $$
(recall that $\gs_{_\el}$ is the value of $s$ when $\S_\ell$ reaches a step with $s\geq 2r_{_{\el}}$ for the first time). If $\S_{\ell}$ is not good, then $\gs_{_\el}=F_{\ell}=0$.
We estimate  $F_\ell  $ rather than $\gs_{_\el}$, as $ F_\ell$ can be bounded more nicely.}

 After $\S_\ell$ reaches a step with $s\geq 2r_{_{\el}}$ for the first time, the probability that $\S_{\ell}$ moves forward is at most $\frac{1}{8(i^2+1)}$ until it ends. Moreover, the bound for the probability holds regardless  of $F_{\ell\,'}$, $\ell' <\ell$.
The same argument as above gives, for a positive integer $k$,
$$
 \pr [
 F_\ell  =k| F_{1}, ..., F_{\ell-1}]
 \leq  \pr \Big[ ~\exists\,  k
 ~\mbox{forward steps among the
 first $t$ steps of $\S^*_\ell$}\, \Big]
  $$
 and, by ${t \choose k} \leq (\frac{et}{k})^k$, $\ln (e/8)\leq -1$
 and $\ln (1+y) \leq y$ for $y>0$,
 $$ \pr [ F_\ell =k| F_{1}, ..., F_{\ell-1}]
 \leq  {t \choose k} \Big( \frac{1}{8(i^2+1)}\Big)^k
 \leq \exp\Big( k \ln \frac{e t}{8k(i^2+1)} \Big)
 \leq e^{-k+1}. $$
 The inequality still holds when $k=0$.
For $h=h_q =\lfloor 2^{-(i+1)}q +
q^{3/4} \rfloor $,
$$ \pr \Big[  F_1=k_1, ..., F_{h} =k_h \Big]
=\prod_{\ell=1}^{h} \pr \Big[ F_\ell=k_{\ell} \Big | F_1=k_1, ..., F_{\ell-1} =k_{\ell-1} \Big]
\leq e^{-(\sum_{\ell=1}^{h} k_{\ell})+h},$$
implies that
$$ \pr \Big[ \sum_{\ell=1}^{h} F_\ell =k \Big]
= \!\! \! \sum_{k_{_{\el}}\geq 0 \atop k_{_1}+\cdots+ k_{_h}=k}
\pr \Big[ F_1=k_1, ..., F_{h} =k_h \Big]
\leq {k+h \choose h} e^{-k+h}.$$
Since ${k+h \choose h} \leq (\frac{e(k+h)}{h})^h$, we have
$$
\pr \Big[ \sum_{\ell=1}^{h} F_\ell =k \Big]
\leq \exp\Big( h\ln \frac{e(k+h)}{h} -k +h \Big)
= \exp\Big( h\ln \frac{(k+h)}{h} -k +2h \Big). $$
For $k\geq 4h-1$,
$h\ln \frac{(k+h)}{h} -k +2h\leq - 4k/5 +3h$ yields that
 $$ \pr \Big[ \sum_{\ell=1}^{h} F_\ell \geq 4h-1 \Big]
 =\sum_{k=4h-1}^{\infty} \pr \Big[ \sum_{\ell=1}^{h} F_\ell =k \Big]
\leq \sum_{k=4h-1}^{\infty} e^{-4k/5+3h} \leq 2e^{-(h-4)/5}. $$

\old{Since there are at most $h$ good $S_\ell$'s by (b) of Corollary \ref{m2},
$$ \pr \Big[ \sum_{\ell=1}^{\infty} F_\ell \geq 4h-1 \Big]=
\pr \Big[ \sum_{\ell=1}^{h} F_\ell \geq 4h-1 \Big]
 \leq 2e^{-(h-4)/5}.
$$}

Finally, for good $\S_\ell$, the number of forward steps in $\S_\ell$
is $$ \frac{\gs_{_\el}-(2r_{_{\el-1}}-2)}{2i^2}
+F_{\ell}= \frac{r_{_{\el}}-r_{_{\el-1}}}{i^2} + \frac{\gs_{_\el}/2-r_{_{\el}}+1}{i^2}+F_{\ell}, $$
while the number backward steps in $\S_\ell$ is,
$$\frac{1}{2} \Big(
\gs_{_\el} -(2r_{_{\el}}-2)+ 2i^2 F_{\ell}\Big)=
\gs_{_\el}/2 -r_{_{\el}}+1+i^2  F_{\ell}. $$
As $\gs_{_\el}/2 -r_{_\el} \leq i^2-1$ by \raf{tl},
$$ |\S_\ell| \chi_{_\el}\leq \frac{r_{_{\ell}}-r_{_{\ell-1}}}{i^2}+1+
F_{\ell}+ i^2+ i^2 F_{\ell}=\frac{r_{_{\ell}}-r_{_{\ell-1}}}{i^2}+
(i^2 +1)(F_{\ell}+1). $$
Therefore,
\old{ as $\S_\ell'$ is good for $\ell'< \ell$ if $\S_\ell$ is good,
and $r_{_{\el}}> r_{_{\el-1}}$ if they are defined,}
$$ \pr\Big[ \sum_{\ell=1}^{h} |\S_{\ell}|\chi_{_\el} \geq \,
\frac{r^{*}}{i^2} +4(i^2+1)h
\Big] \leq \pr\Big[ \sum_{\ell=1}^{h} F_{\ell}\geq 4h-1\Big]
\leq 2e^{-(h-4)/5}\leq 2e^{-q^{3/4}/5+1},
$$
where $r^{*}= \max\{ r_{_{\el}}: \mbox{$S_\ell$ is good} \}$.

Suppose every $S_{\ell}$ is good if defined. Then there are
 $ \lceil \frac{2|J|-(2r^*-2)+8i^2 \log q}{2i^2}\rceil$
 more steps after $u_{_{2r^{\!*}}}$ is corrected, and the number $|\S|$  of all steps in $\S$, or equivalently in (iii) for fixed $i$, is
 $$ \Big\lceil \frac{|J|-r^*+1+4i^2 \log q}{i^2}\Big\rceil
 + \sum_{\ell=1}^{h} |\S_{\ell}|\chi_{_\el}
 \leq \frac{|J|-r^*+1}{i^2} + 4\log q +1 + \sum_{\ell=1}^{h}
 |\S_{\ell}|\chi_{_\el}.
$$
Thus, if $\sum_{\ell=1}^{h} |\S_{\ell}|\chi_{_\el}
< \frac{r^*}{i^2} + 4(i^2+1)h$, then
$$|\S| < \frac{|J|}{i^2}+ 4(i^2+1)h  + \frac{1}{i^2}+4\log q+1
<  \frac{|J|}{i^2}+ 4(i^2+2)\Big(\frac{q}{2^{i+1}}+q^{3/4}\Big). $$

By the contrapositive, if $|\S| \geq \frac{|J|}{i^2}+ 4(i^2+2)\Big(\frac{q}{2^{i+1}}+q^{3/4}\Big)$, then either  there is $\S_\ell$ that is defined but not good or
$$\sum_{\ell=1}^{h} |\S_{\ell}|\chi_{_\el}
\geq \, \frac{r_{_{\el}}}{i^2} + 4(i^2+1)h,$$
which gives
\begin{eqnarray*}  \pr \Big[ |\S| \geq \frac{|J|}{i^2}+ 4(i^2+2)\Big(\frac{q}{2^{i+1}}+q^{3/4}\Big)\Big]
=O( q^{-3}+ e^{-q^{3/4}/5 } ) = O( q^{-3}).\end{eqnarray*}

\mn
\qed

\mn

\mn

\section{Finding Weighted Graphs}\label{sgfp}

In this section, we present a randomized algorithm finding weighted
graphs using additive queries, where an additive query asks the sum
of weights of edges with both ends in a fixed set. The algorithm
uses coin weighing algorithms presented in the previous section.

Let $G=(V,E,w_G)$ be a weighted graph with $w_G(e)\not= 0$ for all
$e\in E$. We just say graphs for weighted graphs. First of all, it
is enough to consider bipartite graphs: For general graphs, one may
consider two disjoint copies $X,Y$ of $V$. The copy of $u\in V $ in
$X$ and the copy of $v\in V $ in $Y$ form an edge if and only if
$uv$ is an edge in $G$, and, of course, the weight is inherited.
Then a query of type $w(A, B):=\sum_{x\in A, y\in B} w(x,y)$,
$A\subset X, B \subset Y$ is a linear combination of four additive
queries in $G$, that is,
\beq \label{bip}  w(A, B)=
w_{_{\!G}} (A\cup B) - w_{_{\!G}}(A  \setminus B)  -
w_{_{\!G}}(B\setminus A) + w_{_{\!G}} (A \cap B). \enq In the rest of
this section, we consider weighted bipartite graphs $G=(X\cup Y, E,
w)$ with $|X|=|Y|=n$ and $|E|\leq m$.
 A query means that one takes two sets $A\subset
X$ and $B\subset Y$ and finds out $w(A,B):=\sum_{a\in A, \in B}
w(a,b)$.

\mn

If  $\bO (m\log n)$ queries are allowed, it is easy to find the
graph using the randomized binary search:

\mn
{\bf Randomized Binary Search for Graph} Suppose $n, m\geq 1$ and a bipartite
graph $G=X\cup Y$ with at most $m$ edges and $|X|,|Y|\leq n$   is given.
Then, take random subsets
$X', Y'$ of $X$ and $Y$, respectively, so that each vertex $x\in X$
($y\in Y$, resp.) in $X'$ ($Y'$, resp.) with probability $1/2$,
independently of all other vertices. If $w(X',Y')\not=0 $, find an
edge there using the deterministic binary search. Otherwise, take
a new random sets $X', Y'$ and do it again. Stop when $( 2\lceil \log n \rceil+5)m$ queries are asked. Output all edges found.

\mn
The deterministic binary search means that divide $X'$ into two parts
$X'_1, X'_2$ with size difference at most $1$. If $w(X'_1, Y') \not=0$ take $X'_1$, otherwise, take $X'_2$. Keep doing this until a vertex
$x$ with $w(x,Y') \not= 0$ is found. Then, find $y\in Y'$ with $w(x,y)\not= 0$ using the same method.

\mn

If there is an edge in $G$,  the
probability of $w(X',Y')\not=0 $ is at least $1/4$.
 It may be shown that
$(2\lceil \log n \rceil+4+o(1))m$ queries are enough to find all
edges in $G$, with high  probability.  We may prove
$( 2\lceil \log n \rceil+5)m$ queries are enough with probability
$1-e^{-\Omega(m)}$, a proof of which is presented in Appendix.

\begin{lemma}\label{rbs2} The randomized binary search finds all edges of $G$ with probability $1-e^{-\Omega(m)}$.
\end{lemma}

\mn

For a better query complexity, a more sophisticated algorithm is needed.
We first present an algorithm finding all edges of $G$ when the maximum
degree of $G$ is small, say at most $m^{0.1}$.
Then another algorithm is introduced to find vertices of
large degree and edges containing them. Concatenating two algorithms,
the following theorem may be shown.

\begin{theorem}\label{gfp}
Let $n,m$ be positive integers with $n^2\geq m\geq 2$ and let $\ga, \gb>0$ be  positive real numbers (not necessarily constants) with $2\ga<\gb$.
Suppose  a bipartite (weighted) graph $G$ is given such that each part of $G$ has at most $n$ vertices
and there are $m$ or less edges in $G$.
 If the weights $w(e)$ of edges satisfy $\ga \leq |w(e)|\leq \gb$, then there is a randomized polynomial time algorithm that asks $O(\frac{m \log (\gb/\ga) \log n}{\log m})$ queries, and finds all edges with probability $1-O(1/{m^{0.02}})$.
\end{theorem}

\mn
Theorem \ref{gfpm} follows from the theorem and \raf{bip}.

\old{\mn {\bf Finding Large Degree Vertices} To improve the query
complexity, we first find large degree vertices and edges containing
them. .}

\mn

For the first algorithm, let $\gd=0.05$ and assume that the maximum degree of $G$ is less than $m^{2\gd}$.
To present the algorithm, construct a random partition $X_1,
..., X_{m^{1/2+2\gd}}$ of $X$ so that each vertex $x\in X$ is
equally likely in $X_j$, $j=1,...,m^{1/2+2\gd}$, independently of
all other vertices. Similarly, construct a random partition $Y_1,
..., Y_{m^{1/2+2\gd}}$ of $Y$.

\begin{lemma}\label{random2} Under the same hypotheses as in Theorem \ref{gfp}, if the maximum degree of $G$ is less than $m^{2\gd}$, then, with probability $1-(1+o(1)){m^{-\gd}}$, the followings hold.

\mn
(a) For each $i=1,...,m^{1/2+2\gd}$,  $|N(X_i) |\leq 2m^{1/2-2\gd}$, where $N (X_i)  :=\{ y\in Y:y\sim x ~\mbox{for some $x\in X_i$}\}$.

\mn (b) For each $i=1,...,m^{1/2+2\gd}$ and  $y\in Y$, $d(y; X_i) \leq 3$, where  $d(y;X_i):=
\{ x\in X_i : x\sim y \}$.

\mn (c) For each $i=1,...,m^{1/2+2\gd}$, the number of vertices $y\in Y$ with $d(y;X_i)\geq 2$ is at most $m^{5\gd}$.

\mn (d) The statements (a)-(c) hold after the roles of $X$ and $Y$ are switched.

\mn (e) Except for  $3m^{1-3\gd}$ edges, every edge is a unique edge between $X_i$ and $Y_j$ for some pair $i,j$.
\end{lemma}

\proof Let $p=m^{-1/2-2\gd}$. Then,  $\pr[x\in X_i] =p$ for all $x$ and $i$. It is enough to show  that (a)-(c) hold with probability $1-o(m^{-\gd})$ and (e) holds  with probability $1-m^{-\gd}$.

For (a), as
$$ E[ |N(X_i)|] = \sum_{y\in Y} \Big(1-\pr[ X_i \cap N(y) =\emptyset]\Big)
=\sum_{y\in Y} \Big( 1-(1-p)^{d(y)}\Big) \leq
\sum_{y\in Y} pd(y) \leq pm = m^{1/2-2\gd}, $$
 the generalized martingale inequality
(Lemma \ref{gm}) with $p=m^{-1/2-2\gd}$, $c_{_x} = d(x)$, $\lambda =
m^{1/2-2\gd}$, and $\rho=m^{-2\gd}/2$, gives that
$$ \Pr\Big[ |N_H(X_i)| \geq 2m^{1/2-2\gd} \Big] \leq 2\exp\Big(
-\frac{m^{1/2-4\gd}}{2} +\frac{m^{-1/2-6\gd}}{8} \sum_{x\in X }
(d(x))^2 e^{m^{-2\gd}d(x)/2}\Big).$$ Since $e^{m^{-2\gd}d(x)/2} \leq
e^{1/2}\leq 2$ and $\sum_{x\in X } (d(x))^2 \leq m^{2\gd} \sum_{x\in X }
d(x) = m^{1+2\gd}$, we have
$$ \Pr\Big[ |N (X_i)| \geq 2m^{1/2-2\gd} \Big] \leq 2\exp\Big(
-\frac{m^{1/2-4\gd}}{4}\Big), $$
and
$$ \Pr\Big[\, \exists\, i~~ {\rm s.t.} ~ ~ |N (X_i)| \geq 2m^{1/2-2\gd} \Big] \leq 2m^{1/2+2\gd}\exp\Big(
-\frac{m^{1/2-4\gd}}{4}\Big)=o(m^{-\gd}). $$

For (b),
$$ \pr [ d(y; X_i) \geq 4] \leq { d(y) \choose 4} p^4
\leq \frac{(pd(y))^4}{24}. $$
Thus, the probability that  there is  a pair $y,j$ such that
$d(y, X_j) \geq 4$ is at most
$$ \sum_{j=1}^{m^{1/2+2\gd}}\sum_{y\in Y}   \frac{(pd(y))^4}{24}
\leq \frac{p^4  m^{1/2+2\gd} m^{6\gd}}{24}\sum_{y\in Y} d(y)
\leq \frac{m^{-2-8\gd} m^{1/2+2\gd} m^{1+6\gd}}{24}=
\frac{1}{24m^{1/2}}=o(m^{-\gd}). $$

\newcommand{\0}{\emptyset}
For (c), suppose the number $Z_i$ of vertices $y\in Y$ with
$d(y, X_i)\geq 2$ is more than $m^{5\gd}$. Then there are distinct
vertices $y_1,..., y_{m^{\gd}}$ in $Y$ with $d(y_j, X_i)\geq 2$,
$j=1,..., m^{\gd}$, such that $N(y_j)\cap N(y_k)=\0$
for all distinct pairs $j,k=1,..., m^{\gd}$.
 This is possible since each fixed $y\in Y$ satisfies
 $N(y)\cap N(y')\not=\0$ for at most $m^{4\gd}-1$ vertices $y'\in Y$.
As $r!\geq (\frac{r}{e})^r$ and $(d(y_j))^2 \leq m^{2\gd}d(y_j)$,
$$ \pr[ Z_i > m^{5\gd}] \leq \frac{1}{m^{\gd}!}
\sum_{y_{_1},..., y_{m^{\gd}}}\prod_{j=1}^{m^{\gd}} p^2{d(y_j) \choose 2}
\leq \Big(\frac{ep^2m^{2\gd}}{2m^{\gd}}\Big)^{m^{\gd}}
\sum_{y,..., y_{m^{\gd}}} \prod_{j=1}^{m^{\gd}} d(y_j)
\leq \Big(\frac{ep^2m^{\gd}m}{2}\Big)^{m^{\gd}}
$$ 
and
$$ \pr[\, \exists \, i ~~{\rm s. t.}~~ Z_i > m^{5\gd}]
\leq m^{1/2+2\gd} \Big(\frac{e}{2m^{3\gd}}\Big)^{m^\gd} =o(m^{-\gd}).$$

For (e), the probability that an edge $e=xy$ is not a unique edge between any pair of $X_i$ and $Y_j$ is
$$ \sum_{i,j=1}^{m^{1/2+2\gd}}
\Pr[ (x,y) \in X_i\times Y_j ] \pr\Big[ \mbox{$\exists\! \!$~\  edge between $X_i$ and $Y_j$ other than $e$}\Big|(x,y) \in X_i\times Y_j
\Big].
$$
Since  the conditional probability is at most
$$
 (d(x)-1)p + (d(y)-1)p + (m-d(x)-d(y)+1)p^2  \leq 2m^{2\gd} p +
mp^2 \leq 3m^{-4\gd}
$$
and $ \sum_{i,j=1}^{m^{1/2+2\gd}}
\Pr[ (x,y) \in X_i\times Y_j ] =1$, the number $W$ of edges that are  not a unique edge between any pair of $X_i$ and $Y_j$ is at most $3m^{1-4\gd}$ in expectation.
Markov inequality implies that $$ \pr[ W \geq 3 m^{1-3\gd}] \leq m^{-\gd}. $$

\mn
\qed

\mn

The next algorithm finds all edges of $G$ when the maximum degree of $G$ is less than $m^{2\gd}$.

\mn
{\bf Algorithm A} (i)
 For each $i$, $i=1,..., m^{1/2+2\gd}$, regarding each
$y\in Y$ as a coin with weight $w_i (y):=w_G (X_i,y)=\sum_{x\in X_i} w_{G} (x,y)$,
apply the coin weighing algorithm in Corollary \ref{cwc} to find all counterfeit coins with parameters $(m, n, \ga, \gb,\eps, \mu) $ replaced by $(2m^{1/2-2\gd}, n, \ga, 3\gb, m^{-1/2+7\gd}, \frac{4}{1-4\gd})$.
Let $N_0 (X_i)$ be the set of all counterfeit coins found, $i=1,..., m^{1/2+2\gd}$. Do the same for $Y_j$ and let $N_0 (Y_j)$ be  the set  of all counterfeit coins found, $j=1,..., m^{1/2+2\gd}$.

\mn (ii) For all pairs $i,j=1,..., m^{1/2+2\gd}$ with $|N_0 (X_i) \cap
Y_j |= |X_i \cap N_0 (Y_j)|=1$, take $y\in N_0 (X_i) \cap Y_j$ and
$x\in X_i \cap N_0 (Y_j)$ and weigh the possible edge $xy$ to obtain
$w_{G} (x,y)$. For each pair $xy$ with $w_{G} (x,y)\not=0$, declare that $xy$ is an edge of $G$

\mn (iii) Find remaining edges one by one by applying the randomized
 binary search using  no more than $(6\lceil \log n \rceil+15) m^{1-3\gd} $ queries.

\mn

\mn

For the collectedness and the query complexity of the algorithm, we prove the following lemma.

\begin{lemma}\label{aa} Under the same hypotheses as in Theorem \ref{gfp}, if the maximum degree of $G$ is less than $m^{2\gd}$, then, with probability $1-(1+o(1)){m^{-\gd}}$, Algorithm A asks $O(\frac{ m \log (\gb/\ga) \log n}{\log m})$ queries to find all edges of $G$.
\end{lemma}

\proof  Suppose (a)-(e) of Lemma \ref{random3} hold. First, we show that
the parameters $(2m^{1/2-2\gd}, n, \ga, 3\gb, m^{-1/2+7\gd})$
satisfy all the requirements in Corollary \ref{cwc}.
If $y$ is counterfeit, then $w_i (y)= w_G(X_i,y) \not= 0$.
 This gives  $y\in N(X_i)$ and hence  the number of counterfeit coins
 is at most $|N(X_i)| \leq  2m^{1/2-2\gd}$ by (a) of Lemma \ref{random2}.
 The number of all coins is $|Y|\leq n$.
If $y\sim x$ for only one $x \in X_i$, then $|w_i(y)|= |w_G(x,y)|\geq \ga$.
Thus, $0<|w_i(y)|< \ga$ implies $d(y;X_i)\geq 2$. The number of such $y\in Y$
is at most $m^{5\gd}=  m^{-1/2+7\gd}\cdot 2m^{1/2-2\gd}$ by (c) of
Lemma \ref{random2}. Since $d(y;X_i)\leq 3$ by (b) of Lemma \ref{random2}, $|w_i(y)|\leq \sum_{x\in X_i} |w(x,y)| \leq 3\gb$. Therefore, the algorithm finds the set $N_{0} (X_i)$ of all counterfeit coins, with probability $1-O(m^{-2})$ for each $X_i$. Similarly, the algorithm finds the set $N_{0} (Y_j)$ of all counterfeit coins, with probability $1-O(m^{-2})$ for each $Y_j$. Since there are $2m^{1/2+2\gd}$ sets $X_i$ and $ Y_j$,
 $N_{0} (X_i)=\{ y\in Y: w_{i} (y)\not=0\} $ and
 $N_{0} (Y_j)=\{ x\in X: w_{j} (x)\not=0\} )$,
 with probability $1-O(1/m)$.

If  $e=xy$ is a unique edge between $X_i$ and $Y_j$, then $|w_i (y)|, |w_j (x)|\geq \ga $, especially, $y\in N_0( X_i)$ and $x\in N_0 (Y_j)$. Moreover, as there is no other edge between $X_i$ and $Y_j$, $N_0 (X_i) \cap
Y_j =\{y\}$ and $X_i \cap N_0 (Y_j)=\{ x\}$. Thus, the algorithm finds the edge $e=xy$ in (ii).
By (e) of Lemma \ref{random2}, at most $3m^{1-3\gd}$ edges remain unfound in (ii).
All the remaining edges can be found in (iii) with probability $1-e^{-\Omega(m^{1-3\gd})}$ by Lemma \ref{rbs}.

For the query complexity, in (i),  $O(\frac{ m^{1/2-2\gd} \log (\gb/\ga) \log n}{\log m})$ queries are enough for each $X_i$ or $Y_j$. As there are $2m^{1/2+2\gd}$ such sets, $O(\frac{ m \log (\gb/\ga) \log n}{\log m})$ queries are enough in (i). In (ii), if $|N_0 (X_i) \cap
Y_j |= |X_i \cap N_0 (Y_j)|=1$, then there is at least one edge between $X_i$ and $Y_j$.
As there are at most $m$ such pairs $X_i, Y_j$, $m$ queries are enough in (ii).
Since $o(\frac{m\log n}{\log m})$ queries are asked in (iii), the query complexity of the algorithm is
$O(\frac{ m \log (\gb/\ga) \log n}{\log m})$.

\mn\qed

\mn

For general graphs, select each vertex of $Y$ with probability
$m^{-\gd}$, $\gd=0.05$,  independently of all other vertices. Let $G_1$ be the induced graph  on $X$ and the selected
vertices of $Y$.
\begin{lemma}\label{random3} If $G$ has at most $m$ edges, then the followings hold with probability $1-O(m^{-\gd/2})$.

\mn
(a) The number of edges in $G_1$ is at most $m^{1-\gd/2}$.

\mn (b) If $d_{G_1} (x) \geq  m^\gd/2$, then $
d_G (x) \leq 2m^\gd d_{G_1} (x)\leq 3d_{G} (x). $

\mn (c) If $d_{G} (x) \geq  m^{2\gd}$, then $d_{G_1} (x) \geq  m^\gd/2$.

\end{lemma}

\proof As each edge in $G_1$ with probability $m^{-\gd}$,
the expected number of edges in $G_1$ is at most $m^{1-\gd}$.
Markov Inequality gives
$$ \Pr[ \mbox{the number of edges in $G_1$} \geq m^{1-\gd/2} ] \leq
m^{-\gd/2}. $$
For the degree $d_{G_1} (x)$ of $x$ in $G_1$, as $E[d_{G_1}(x)] = m^{-\gd} d_G(x)$, Lemma \ref{gm} with $c_y =1$ if $y\sim x$ and $c_y =0$
otherwise, $\lambda =\frac{m^{\gd} }{4}$, $\rho=1/2$ gives, for $x\in X$ with $d_G(x) < \frac{m^{2\gd}}{4}$,
$$ \Pr \Big[ |d_{G_1} (x)-m^{-\gd}d_G (x)| \geq   \frac{ m^{\gd}}{4}\Big] \leq
2 \exp \Big(-\frac{m^{\gd}}{8}+ \frac{e^{1/2} m^{-\gd}
d_G (x)}{8} \Big) \leq  2\exp \Big(-\frac{m^{\gd} }{16}\Big).
$$
In particular, if $d_G(x) < \frac{m^{2\gd}}{4}$, then
$d_{G_1} (x)-m^{-\gd}d_G (x) < m^\gd/4$, or equivalently,
$d_{G_1} (x) < m^\gd/4+m^{-\gd}d_G (x)< m^\gd/2$, with probability $1-e^{-\Omega(m^\gd)}$.

For (b),  it is now enough to show that $
d_G (x) \leq 2m^\gd d_{G_1} (x)\leq 3d_{G} (x)$ when  $d_G(x) \geq  \frac{m^{2\gd}}{4}$,  say, with probability $1-e^{-\Omega(m^{\gd})}$. Lemma \ref{gm} with $c_y =1$ if $y\sim x$ and $c_y =0$
otherwise, $\lambda =\frac{m^{-\gd} d_G(x)}{2}$, $\rho=1/3$  also gives
$$ \Pr \Big[ |d_{G_1} (x)-m^{-\gd}d_G (x)| \geq   \frac{ m^{-\gd}d_G (x)}{2}\Big] \leq
2 \exp \Big(-\frac{m^{-\gd} d_G(x)}{6}+ \frac{e^{1/3} m^{-\gd}
d_G (x)}{18} \Big) \leq  2\exp \Big(-\frac{m^{-\gd} d_G (x)}{12}\Big),
$$
for $e^{1/3} \leq 3/2$. If $d_G (x)\geq m^{2\gd}/4$, we have
$ |2m^{\gd} d_{G_1}(x) -2d_{G}(x)| \leq d_G (x)$, or equivalently,
 $ d_G (x)\leq 2m^{\gd} d_{G_1}(x) \leq 3d_{G}(x) $,
 with probability $1-e^{-\Omega(m^{\gd})}$.
Moreover, if
 $d_G (x)\geq m^{2\gd}$, then $2m^{\gd} d_{G_1}(x) \geq d_G (x) \geq m^{2\gd}$.  That is, $d_{G_1}(x)\geq m^{\gd}/2$, which shows (c).

\mn
\qed

\bn

\mn
{\bf Algorithm B}
 (i)  Apply the randomized
binary search to find edges of $G_1$ one by one, using $(2\lceil \log n \rceil +5)m^{1-\gd/2}$ queries. Let $G_2$ be the graph on $X\cup Y$ consisting of all edges found.

\mn  (ii) For each vertex  $x\in X$ with  $ d_{G_2} (x) \geq m^\gd /2$, regard each $y\in Y$ as a coin with weight $w_x (y):=w_G (x,y)$ and apply the  coin weighing algorithm in Corollary \ref{cwc} with parameters $(m, n, \ga, \gb,\eps, \mu) $ replaced by $(2m^{\gd} d_{G_2} (x), n, \ga, \gb,0, 1/\gd)$. The vertices $x\in X$  with   $ d_{G_2} (x) \geq m^\gd /2$ are called vertices of large degree.

\mn (iii)
Output vertices of large degree and all edges found.

\mn

\mn

Algorithm B has the following property.

\begin{lemma}\label{ab}  Under the same hypotheses as in Theorem \ref{gfp}, with probability $1-O(m^{-\gd/2})$,
Algorithm B uses $O(\frac{m \log (\gb/\ga) \log n}{\log m})$ queries to find all vertices $x\in X$ with $d_G(x) \geq m^{2\gd}$ and all edges containing them.
\end{lemma}

\proof Suppose (a) and (b) of Lemma \ref{random3} hold. Then
Lemma \ref{rbs2} yields $G_2=G_1$ with probability
$1-e^{\Omega(m^{1-\gd/2})}$. We assume
that $G_1=G_2$ in the rest of the proof.

In (ii), note that the number of counterfeit coins for
$x$ is $d_{G} (x)$, which is at most $2m^{\gd } d_{G_2} (x) $
for all $d_{G_1} (x)=
d_{G_2} (x)\geq m^{\gd}/2$ by (b) of Lemma \ref{random3}.
Thus, the algorithm in Corollary \ref{cwc} finds $N_G (x)$
for each $x\in X$ satisfying $d_{G_2}(x) \geq m^{\gd}/2$,
with probability $1- O(1/(2m^{\gd } d_{G_2} (x))^{1/\gd})=1- O(1/m^2)$.
As $d_G(x)\geq d_{G_2}(x) $, there are at most $2m^{1-\gd} $
vertices $x \in X$ with $d_{G_2}(x) \geq m^{\gd}/2$ and the algorithm
finds   $N_G (x)$ for all such vertices $x\in X$, with probability $1-O(1/m)$.
In particular, if $d_G(x) \geq m^{2\gd}$, then $d_{G_2}(x) \geq
m^{\gd}/2$ by (c) of Lemma \ref{random3} and hence $N_G(x)$ are found.

For the query complexity, $(2\lceil \log n \rceil +5)m^{1-\gd/2}$ queries are asked in (i). In (ii), $O(\frac{ m^{\gd}d_{G_2} (x) \log (\gb/\ga) \log n}{\log m})$ queries are asked for each $x\in X$ with $d_{G_2} (x) \geq m^{\gd}/2$. On the other hand,  $d_{G_2} (x) \geq m^{\gd}/2$ implies  $2m^{\gd}d_{G_2} (x)\leq 3d_{G} (x)$ by (b) of Lemma \ref{random3}. Thus,
 $$\sum_{x: d_{G_2} (x) \geq m^{\gd}/2} m^{\gd}d_{G_2} (x)
  \leq \frac{3}{2} \sum_{x\in X} d_{G} (x)= \frac{3m}{2} $$
 gives that $O(\frac{m \log (\gb/\ga) \log n}{\log m})$ queries are asked in (ii).

 \mn
 \qed

\bn

To find all vertices $v$ in $G$ with $d_G (v) \geq m^{2\gd}$, one may
apply Algorithm B twice, one as it is and the other after exchanging roles
of $X$ and $Y$. Then, after removing all vertices found (and all edges
containing any of them),  we apply Algorithm A.
Lemmas \ref{aa} and \ref{ab} imply that

\begin{corollary}
Under the same hypotheses as in Theorem \ref{gfp}, there is a polynomial time
randomized algorithm asking $O(\frac{ m \log (\gb/\ga) \log n}{\log m})$ queries
 to find all edges of $G$, with probability $1-O(1/m^{0.02})$.
\end{corollary}

If  the algorithm in the corollary is forced to stop when
it asks $\frac{ \eta m \log (\gb/\ga) \log n}{\log m}$ queries,
for  the constant $\eta$   in the $O(\frac{ m \log (\gb/\ga) \log n}{\log m})$
term,
the desired algorithm  in Theorem \ref{gfp} may be obtained.

\mn

\mn
\section{Concluding Remarks}

In this paper, we presented a polynomial time randomized algorithm that
uses $O(\frac{ m \log (\gb/\ga) \log n}{\log m})$ queries, when there are
at most $m$ counterfeit coins and the weights $w(c)$ of all
 counterfeit coins
satisfy $\ga\leq |w(c)|\leq \gb$. This plays a key role to find a hidden
weighted graph $G$    satisfying similar  conditions.  Though there is a non-adaptive algorithm to find all counterfeit coins using $O(\frac{ m \log n}{\log m})$ queries \cite{BM11_TCS}, it is not a polynomial time algorithm.
An obvious question is if  there is a polynomial time algorithm to find all
counterfeit coins using $O(\frac{ m \log n}{\log m})$ queries when there is
no restriction on the wights.

 The algorithm we presented was a randomized algorithm that uses the optimal
 number of queries up to a constant factor. On the other hand, the best
 deterministic algorithm uses $\Theta (\frac{m\log n}{\log m} + m\log\log m)$
 (see \cite{BM10_MFCS}), it would be good to implement a deterministic
 polynomial time algorithm that uses $O(\frac{ m \log n}{\log m})$
 queries even when the weights of counterfeit coins are positive real numbers.

\bibliographystyle{plain}
\bibliography{biblio}  

\begin{thebibliography}{10}

\bibitem{AA05}
N.~Alon and V.~Asodi.
\newblock Learning a hidden subgraph.
\newblock {\em SIAM Journal on Discrete Mathematics}, 18(4):697--712, 2005.

\bibitem{ABKRS04}
N.~Alon, R.~Beigel, S.~Kasif, S.~Rudich, and B.~Sudakov.
\newblock Learning a hidden matching.
\newblock {\em SIAM Journal on Computing}, 33(2):487--501, 2004.

\bibitem{AC04}
D.~Angluin and J.~Chen.
\newblock Learning a hidden graph using $\mathcal{O}(\log{n})$ queries per
  edge.
\newblock In {\em Proceedings of the 17th Annual Conference on Learning Theory
  (COLT 2004)}, pages 210--223, Banff, Canada, 2004.

\bibitem{AC06}
D.~Angluin and J.~Chen.
\newblock Learning a hidden hypergraph.
\newblock {\em Journal of Machine Learning Research}, 7:2215--2236, 2006.

\bibitem{BAFK01}
R.~Beigel, M.~S. Apaydin, L.~Fortnow, and S.~Kasif.
\newblock An optimal procedure for gap closing in whole genome shotgun
  sequencing.
\newblock In {\em Proceedings of the Fifth Annual International Conference on
  Computational Molecular Biology (RECOMB 2001)}, pages 22--30, 2001.

\bibitem{BGK05}
M.~Bouvel, V.~Grebinski, and G.~Kucherov.
\newblock Combinatorial search on graphs motivated by bioinformatics
  applications: A brief survey.
\newblock In {\em the 31st International Workshop on Graph-Theoretic Concepts
  in Computer Science (WG 2005)}, pages 16--27, 2005.

\bibitem{Bshouty09}
N.~H. Bshouty.
\newblock Optimal algorithms for the coin weighing problem with a spring scale.
\newblock In {\em Proceedings of the 22nd Annual Conference on Learning Theory
  (COLT 2009)}, Montreal, Canada, 2009.

\bibitem{BM10_STACS}
N.~H. Bshouty and H.~Mazzawi.
\newblock Optimal query complexity for reconstructing hypergraphs.
\newblock In {\em Proceedings of the 27th International Symposium on
  Theoretical Aspects of Computer Science (STACS 2010)}, pages 143--154, Nancy,
  France, 2010.

\bibitem{BM10_MFCS}
N.~H. Bshouty and H.~Mazzawi.
\newblock Toward a deterministic polynomial time algorithm with optimal
  additive query complexity.
\newblock In {\em Proceedings of the 35th International Symposium on
  Mathematical Foundations of Computer Science (MFCS 2010)}, pages 221--232,
  Brno, Czech Republic, 2010.

\bibitem{BM11_SODA}
N.~H. Bshouty and H.~Mazzawi.
\newblock On parity check $(0,1)$-matrix over $\mathbb{Z}_p$.
\newblock In {\em Proceedings of the ACM-SIAM Symposium on Discrete Algorithms
  (SODA 2011)}, pages 1383--1394, San Francisco, USA, 2011.

\bibitem{BM11_TCS}
N.~H. Bshouty and H.~Mazzawi.
\newblock Reconstructing weighted graphs with minimal query complexity.
\newblock {\em Theoretical Computer Science}, 412(19):1782--1790, 2011.

\bibitem{Cantor64}
D.~G. Cantor.
\newblock Determining a set from the cardinalities of its intersections with
  other sets.
\newblock {\em Canadian Journal of Mathematics}, 16:94--97, 1964.

\bibitem{CM66}
D.~G. Cantor and W.~H. Mills.
\newblock Determination of a subset from certain combinatorial properties.
\newblock {\em Canadian Journal of Mathematics}, 18:42--48, 1966.

\bibitem{Capetanakis79b}
J.~Capetanakis.
\newblock Generalized {TDMA}: The multi-accessing tree protocol.
\newblock {\em IEEE Transactions on Communications}, 27(10):1476--1484, 1979.

\bibitem{Capetanakis79a}
J.~Capetanakis.
\newblock Tree algorithms for packet broadcast channels.
\newblock {\em IEEE Transactions on Information Theory}, 25(5):505--515, 1979.

\bibitem{CJK11_JCSS}
S.~S. Choi, K.~Jung, and J.~H. Kim.
\newblock Almost tight upper bound for finding {F}ourier coefficients of
  $k$-bounded pseudo-{B}oolean functions.
\newblock {\em Journal of Computer and System Sciences}, 77(6):1039--1053,
  2011.

\bibitem{CK12}
S.~S. Choi and J.~H. Kim.
\newblock Randomized polynomial time algorithms for finding weighted graphs
  with optimal additive query complexity.
\newblock {\em {\em submitted}}.

\bibitem{CK10_AI}
S.~S. Choi and J.~H. Kim.
\newblock Optimal query complexity bounds for finding graphs.
\newblock {\em Artificial Intelligence}, 174(9--10):551--569, 2010.

\bibitem{CK11_AI}
S.~S. Choi and J.~H. Kim.
\newblock Sample complexity for linkage learning.
\newblock Submitted, 2011.

\bibitem{DH93}
D.~Du and F.~K. Hwang.
\newblock Combinatorial group testing and its application.
\newblock In {\em V. 3 of Series on applied mathematics}, chapter~10. World
  Science, 1993.

\bibitem{ER63}
P.~Erd\H{o}s and A.~R\'{e}nyi.
\newblock On two problems of information theory.
\newblock {\em Publications of the Mathematical Institute of the Hungarian
  Academy of Sciences}, 8:241--254, 1963.

\bibitem{Fine60}
N.~J. Fine.
\newblock Solution of problem {E} 1399.
\newblock {\em American Mathematical Monthly}, 67(7):697--698, 1960.

\bibitem{Grebinski98}
V.~Grebinski.
\newblock On the power of additive combinatorial search model.
\newblock In {\em Proceedings of the 4th Annual International Conference on
  Computing and Combinatorics (COCOON 1998)}, pages 194--203, Taipei, Taiwan,
  1998.

\bibitem{GK97}
V.~Grebinski and G.~Kucherov.
\newblock Optimal query bounds for reconstructing a {H}amiltonian cycle in
  complete graphs.
\newblock In {\em the Fifth Israel Symposium on the Theory of Computing Systems
  (ISTCS 1997)}, pages 166--173, 1997.

\bibitem{GK98}
V.~Grebinski and G.~Kucherov.
\newblock Reconstructing a {H}amiltonian cycle by querying the graph:
  Application to {DNA} physical mapping.
\newblock {\em Discrete Applied Mathematics}, 88:147--165, 1998.

\bibitem{GK00}
V.~Grebinski and G.~Kucherov.
\newblock Optimal reconstruction of graphs under the additive model.
\newblock {\em Algorithmica}, 28:104--124, 2000.

\bibitem{Hein89}
J.~J. Hein.
\newblock An optimal algorithm to reconstruct trees from additive distance
  data.
\newblock {\em Bulletin of Mathematical Biology}, 51(5):597--603, 1989.

\bibitem{Kim95}
J.~H. Kim.
\newblock The {R}amsey number ${R}(3,t)$ has order of magnitude $t^2/\log{t}$.
\newblock {\em Random Structures and Algorithms}, 7(3):173--208, 1995.

\bibitem{KZZ03}
V.~King, L.~Zhang, and Y.~Zhou.
\newblock On the complexity of distance-based evolutionary tree reconstruction.
\newblock In {\em Proceedings of the 14th Annual ACM-SIAM Symposium on Discrete
  Algorithms (SODA 2003)}, pages 444--453, 2003.

\bibitem{LV94}
M.~Li and P.~M.~B. Vit{\'{a}}nyi.
\newblock Kolmogorov complexity arguments in combinatorics.
\newblock {\em J. Comb. Theory Series A}, 66(2):226--236, 1994.

\bibitem{Lindstrom64}
B.~Lindstr{\"{o}}m.
\newblock On a combinatory detection problem {I}.
\newblock {\em Publications of the Mathematical Institute of the Hungarian
  Academy of Sciences}, 9:195--207, 1964.

\bibitem{Lindstrom65}
B.~Lindstr{\"{o}}m.
\newblock On a combinatorial problem in number theory.
\newblock {\em Canadian Mathematical Bulletin}, 8(4):477--490, 1965.

\bibitem{Lindstrom71}
B.~Lindstr{\"{o}}m.
\newblock On {M}{\"{o}}bius functions and a problem in combinatorial number
  theory.
\newblock {\em Canadian Mathematical Bulletin}, 14(4):513--516, 1971.

\bibitem{Lindstrom75}
B.~Lindstr{\"{o}}m.
\newblock Determining subsets by unramified experiments.
\newblock In J.~N. Srivastava, editor, {\em A Survey of Statistical Designs and
  Linear Models}, pages 407--418. North Holland, 1975.

\bibitem{M81}
J.~L. Massey.
\newblock Collision-resolution algorithms and random-access communications.
\newblock In G.~Longo, editor, {\em Multi-user communications systems, CISM
  Courses and Lecture Notes No. 265}, pages 73--137. Springer, Wien and New
  York, 1981.

\bibitem{Mazzawi10}
H.~Mazzawi.
\newblock Optimally reconstructing weighted graphs using queries.
\newblock In {\em Proceedings of the ACM-SIAM Symposium on Discrete Algorithms
  (SODA 2010)}, pages 608--615, Austin, USA, 2010.

\bibitem{McDiarmid89}
C.~McDiarmid.
\newblock On the method of bounded differences.
\newblock In J.~Siemons, editor, {\em Surveys in Combinatorics}, London
  Mathematical Society Lecture Note Series 141, pages 148--188. Cambridge
  University Press, 1989.

\bibitem{Moser70}
L.~Moser.
\newblock The second moment method in combinatorial analysis.
\newblock In {\em Combinatorial Structures and Their Applications. Proceedings
  of the Calgary International Conference on Combinatorial Structures and Their
  Applications held at the University of Calgary. June 1969}, pages 283--384.
  Gordon and Breach, New York, 1970.

\bibitem{RS07}
L.~Reyzin and N.~Srivastava.
\newblock Learning and verifying graphs using queries with a focus on edge
  counting.
\newblock In {\em Proceedings of the 18th International Conference on
  Algorithmic Learning Theory (ALT 2007)}, pages 285--297, Sendai, Japan, 2007.

\bibitem{RS07a}
L.~Reyzin and N.~Srivastava.
\newblock On the longest path algorithm for reconstructing trees from distance
  matrices.
\newblock {\em Information Processing Letters}, 101(3):98--100, 2007.

\bibitem{Shapiro60}
H.~S. Shapiro.
\newblock Problem {E} 1399.
\newblock {\em American Mathematical Monthly}, 67(1):82, 1960.

\bibitem{SS63}
S.~S{\"{o}}derberg and H.~S. Shapiro.
\newblock A combinatory detection problem.
\newblock {\em American Mathematical Monthly}, 70:1066--1070, 1963.

\bibitem{TRKKS99}
H.~Tettelin, D.~Radune, S.~Kasif, H.~Khouri, and S.~L. Salzberg.
\newblock Optimized multiplex {PCR}: Efficiently closing a whole-genome shotgun
  sequencing project.
\newblock {\em Genomics}, 62:500--507, 1999.

\bibitem{TM78}
B.~Tsybakov and V.~Mikhailov.
\newblock Free synchronous packet access in a broadcast channel with feedback.
\newblock {\em Problemy Peredachi Informassi}, 14(4):259--280, 1978.

\bibitem{UTW00}
R.~Uehara, K.~Tsuchida, and I.~Wegener.
\newblock Identification of partial disjunction, parity, and threshold
  functions.
\newblock {\em Theoretical Computer Science}, 210(1--2):131--147, 2000.

\end{thebibliography}

\mn

\mn

\mn

\mn

\mn
{\LARGE \bf Appendix: Proofs of Lemmas \ref{rbs}, \ref{random} and \ref{rbs2} }

\mn

\mn

In this appendix, we prove Lemmas \ref{rbs}, \ref{random} and \ref{rbs2}.

\mn

\mn
{\bf Lemma \ref{rbs}} {\em The randomized binary search finds all counterfeit coins with probability $1-e^{-\Omega(m)}$.
}

\proof If there is a counterfeit coin $c$, conditioned on $A''=A'\setminus \{c\}$, $A'$ can be one of $A''$ and $A''\cup \{ c\}$, each with probability $1/2$. Since $w(A''\cup \{ c\})=w(A'')+ w(c) \not= w(A'')$, the probability of $ w(A')\not=0$ is at least   $1/2$.

Let $Z_i$ be the number of random trials when the $i^{\rm th}$ counterfeit coin is found. Then, for $a=\ln(4/3)$,
$$ E[e^{aZ_i}| Z_1,...,Z_{i-1}]
= \sum_{k=1}^{\infty} e^{ak} (1-p_{_i})^{k-1} p_{_i}
= \frac{p_{_i}}{1-p_{_i}} \frac{ (1-p_{_i})e^a}{1-(1-p_{_i})e^a}
=\frac{p_{_i}e^a  }{1-(1-p_{_i})e^a}\leq \frac{e^a  }{2-e^a},$$
for $p_{_i} := \pr [ w(A') \not=0| Z_1,...,Z_{i-1}] \geq 1/2$.
As
$$ E\Big[ e^{a \sum_{i=1}^{\ell} Z_i}\Big]
=  E\Big[ E\Big[e^{a \sum_{i=1}^{\ell} Z_i}\Big|Z_1,...,Z_{\ell-1}\Big]\Big]
=E\Big[e^{a \sum_{i=1}^{\ell-1} Z_i} E\Big[e^{aZ_\ell}\Big| Z_1,...,Z_{\ell-1}\Big]\Big]
\leq \Big(\frac{e^a  }{2-e^a}\Big) E\Big[e^{a \sum_{i=1}^{\ell-1} Z_i} \Big],
$$
for all $\ell=1,..., m$, we have
$$ \pr [ \sum_{i=1}^{m^*} Z_i \geq 3m ]
\leq \pr [ e^{a \sum_{i=1}^{m^*} Z_i} \geq e^{3am}]
\leq E[ e^{a (\sum_{i=1}^{m^*} Z_i-3m)}]
\leq  \Big(\frac{e^{-2a} }{2-e^a}\Big)^{m}\leq \Big(\frac{27}{32}\Big)^m=e^{-\Omega(m)}, $$
where $m^* \leq m$ is the number of counterfeit coins,

\mn
\qed

\mn
{\bf Lemma \ref{random}} {\em
Suppose a set $A$ of   $n$ or less coins are given, and the number of counterfeit coins in $A$ is at most $q\geq 2$. If the weights $w(c)$ of all but at most $q/2$ counterfeit coins $c$ satisfy $|w(c)|\geq \ga$.
 Then,
with
probability $1-O(\frac{1}{q})$, we have the followings.

\mn
(a) There are at most $\frac{5q}{6}$  counterfeit coins $c$ that satisfy $|w(c)|< \ga$  (not exclusively) or belong to  a set $A_{0,j}$ containing more than one counterfeit coin, $j=1,...,2^{\lmq}$.

\mn
(b) For each $i=1,..., \lceil 2\log q \rceil-1$, $A_{i,j}$ contains
at most   $\frac{i+2\log q}{i}$ counterfeit coins.

\mn (c) For each  $i=1,..., \lceil 2\log q \rceil-1$, there are at most  $2^{-(i+1)}q +
q^{3/4}$ sets $A_{i,j}$ that
 contain more than one counterfeit coin.

\mn  (d) For $i\geq \lceil 2\log q\rceil-1$,   each $A_{i,j}$ contains
one or less counterfeit coin.

\mn (e) Each  $A_{\lceil 3\log n \rceil,j}$ contains at most one coin.
}

\proof  (a) For  any counterfeit coin $c$,
the probability that $c$ belongs to a set $A_{0,j}$ containing another counterfeit coin is at most
$1- (1-2^{-\q})^{q-1} \leq 1-1/e.$
Thus,  the  number
of such counterfeit coins $c$ with $|w(c) |\geq \ga$  is at most $(1-1/e)(q-q_{_1})$ in expectation, where $q_{_1}$ is the number of counterfeit coins $c$ with $|w(c)|< \ga$.
As $q-q_{_1} \geq q/2$ and the number depends only on where counterfeit coins are in, we may apply the Azuma-Hoeffding martingale inequality (Lemma \ref{mar}) with $c_{_\ell}=2$ and $\sum_{\ell} c_{_\ell}^2 \leq 4q$ to deduce that the number of such counterfeit coins $c$ with  $|w(c) |\geq \ga$ is at most $2(q-q_{_1})/3$, with probability $1-e^{-\Omega(q)}$. Thus, with probability $1-e^{-\Omega(q)}$, there are at most $$2(q-q_{_1})/3+q_{_1}=2q/3+q_{_1}/3\leq 2q/3+q/6=5q/6 $$ counterfeit coins $c$ that satisfy $|w(c)|< \ga$ or belong to  a set $A_{0,j}$ containing another counterfeit coin.

 \mn
  (b) For each set $A_{i,j}$, the probability that $A_{i,j}$ contains  $k_{_i}:= \lceil
\frac{i+2\log q}{i} \rceil $ or more   counterfeit coins are bounded from above by
$$ {q \choose k_{_i}} 2^{-k_{_i}(\lmq+i)} \leq
2^{-k_{_i}(\lmq+i)}  q^{k_{_i}}.
$$
Thus, for each $i=1,..., \lceil 2\log q \rceil-1$, the probability that
$A_{i,j}$ contains $k_{_i}$ or more counterfeit coins is at most
$$ 2^{\lmq+i} 2^{-k_{_i}(\lmq+i)}  q^{k_{_i}}
= 2^{-(k_{_i}-1)(\lmq+i)} q^{k_{_i}}\leq 2^{-(k_{_i}-1)i} q
 \leq \frac{1}{q}.
$$

\mn
(c)  The probability that $A_{i,j}$ contains two or more counterfeit coins is at most
$$ {q \choose 2} 2^{-2(\q+i)} \leq  2^{-2(\q+i)-1}q^2, $$
and, for each $i$, the expected number of $A_{i,j}$ containing two or more counterfeit coins is at most
$$ 2^{\q+i} 2^{-2(\q+i)-1}q^2 \leq 2^{-(\q+i)-1} q^2= 2^{-(i+1)}q.
$$
Counting coordinates corresponding to counterfeit coins only, we apply the Azuma-Hoeffding martingale inequality (Lemma \ref{mar}) with $c_{_\ell}=1$ to conclude that, for each $i$,  the  number of $A_{i,j}$  containing two or more counterfeit coins is at most $2^{-(i+1)}q+q^{3/4}$, with probability $1-e^{-\Omega(q^{1/2})}$.

\mn
(d) For $i=\lceil 2\log q \rceil-1$, by the same estimation as in (c), the probability that $A_
{i, j} $ contains two or more counterfeit coins for some $j$ is at most $ 2^{-(i+1)}q\leq 1/q$.  If each $A_
{\lceil 2\log q \rceil-1, j} $ contains at most one counterfeit coin, then so does  each $A_{i,j}$  with $i\geq  \lceil 2\log q \rceil$,  for  every set $A_{i,j}$ with $i\geq \lceil 2\log q \rceil$ is a subset of some $A_{\lceil 2\log q \rceil-1, \ell}$.

\mn
(e)  The statement follows since  each set $A_{i,j}$ with $i\geq \lceil 2\log q \rceil-1$ is deterministically divide into two sets with size difference at most $1$.

\mn
\qed

\mn
{\bf Lemma \ref{rbs2}} {\em The randomized binary search finds all edges of $G$ with probability $1-e^{-\Omega(m)}$.
}

\proof If $G$ has at least one edge, say $e=xy$,
conditioned $X'':=X'\setminus \{x\}$
and $Y''= Y'\setminus \{y\}$, $(X', Y')$ can be one of $(X'', Y'')$,
$(X''\cup \{x\}, Y'')$, $(X'', Y''\cup \{y\})$, and $(X''\cup \{x\}, Y''\cup \{y\})$, each with probability $1/4$. If all three weights $w(X'', Y'')$,
$w(X''\cup \{x\}, Y'')$, $w(X'', Y''\cup \{y\})$ are $0$, then
$w(X''\cup \{x\}, Y''\cup \{y\})=w(e)\not=0$. In other words, at least one of the four weights is non-zero. This yields that  the
probability of $w(X',Y')\not=0 $ is at least $1/4$.

Let $Z_i$ be the number of random trials when the $i^{\rm th}$ edge is found. Then, for $a=\ln(13/12)$,
$$ E[e^{aZ_i}| Z_1,...,Z_{i-1}]
= \sum_{k=1}^{\infty} e^{ak} (1-p_{_i})^{k-1} p_{_i}
= \frac{p_{_i}}{1-p_{_i}} \frac{ (1-p_{_i})e^a}{1-(1-p_{_i})e^a}
=\frac{p_{_i}e^a  }{1-(1-p_{_i})e^a}\leq \frac{e^a  }{4-3e^a},$$
for $p_{_i} := \pr [ w(X',Y') \not=0| Z_1,...,Z_{i-1}] \geq 1/4$.
Thus,
$$ \pr [ \sum_{i=1}^{m^*} Z_i \geq 5m ]
\leq \pr [ e^{a \sum_{i=1}^{m^*} Z_i} \geq e^{5am}]
\leq E[ e^{a (\sum_{i=1}^{m^*} Z_i-5m)}]
\leq  \Big(\frac{e^{-4a} }{4-3e^a}\Big)^{m}\leq (0.97)^m=e^{-\Omega(m)}, $$
where $m^* \leq m$ is the number of edges in $G$, as
$$ E\Big[ e^{a \sum_{i=1}^{\ell} Z_i}\Big]
=  E\Big[ E\Big[e^{a \sum_{i=1}^{\ell} Z_i}\Big|Z_1,...,Z_{\ell-1}\Big]\Big]
=E\Big[e^{a \sum_{i=1}^{\ell-1} Z_i} E\Big[e^{aZ_\ell}\Big| Z_1,...,Z_{\ell-1}\Big]\Big]
\leq \Big(\frac{e^a  }{4-3e^a}\Big) E\Big[e^{a \sum_{i=1}^{\ell-1} Z_i} \Big],
$$
for all $\ell=1,..., m^*$.

\mn
\qed

\end{document}